\documentclass[11pt]{article}

\usepackage{amsmath}
\usepackage{amsfonts}
\usepackage{amssymb}
\usepackage{amsthm}
\usepackage{graphicx}
\usepackage{cite}
\usepackage{caption}
\usepackage{url}
\usepackage{color}

\usepackage{algorithm}
\usepackage{algpseudocode}
\usepackage{mathtools}
\usepackage{subcaption}
\usepackage{multirow}
\usepackage{bbm}
\usepackage[table]{xcolor}
\usepackage{xr}

\makeatletter
\newcommand*{\addFileDependency}[1]{
  \typeout{(#1)}
  \@addtofilelist{#1}
  \IfFileExists{#1}{}{\typeout{No file #1.}}
}
\makeatother

\newcommand*{\myexternaldocument}[1]{%
    \externaldocument{#1}%
    \addFileDependency{#1.tex}%
    \addFileDependency{#1.aux}%
}

\myexternaldocument{ParallelIFGF_Supplementary_Materials}

\newcommand\Tstrut{\rule{0ex}{2.7ex}}         
\newcommand\Bstrut{\rule[-1.3ex]{0pt}{0pt}}   
\newcommand\Strut{\rule{0ex}{0ex}}         
\algblock{ParFor}{EndParFor}
\algnewcommand\algorithmicparfor{\textbf{parallel for}}
\algnewcommand\algorithmicpardo{\textbf{do}}
\algnewcommand\algorithmicendparfor{\textbf{end\ parallel for}}
\algrenewtext{ParFor}[1]{\algorithmicparfor\ #1\ \algorithmicpardo}
\algrenewtext{EndParFor}{\algorithmicendparfor}

\definecolor{MyGreen}{rgb}{0,0.6,0}

\title{Massively Parallelized\\ Interpolated Factored Green Function Method} \author{Christoph Bauinger$^*$ and Oscar
  P. Bruno\footnote{Computing and Mathematical Sciences, Caltech, Pasadena, CA 91125, USA}}
\begin{document}
\date{}
\maketitle
\begin{abstract}
  This paper presents the first parallel implementation of the novel
  ``Interpolated Factored Green Function'' (IFGF) method introduced
  recently for the accelerated evaluation of discrete integral
  operators arising in wave scattering and other areas (Bauinger and
  Bruno, Jour. Computat. Phys., 2021). On the basis of the
  hierarchical IFGF interpolation strategy, the proposed (hybrid
  MPI-OpenMP) parallel implementation results in highly efficient data
  communication, and it exhibits in practice excellent parallel
  scaling up to large numbers of cores---without any hard limitations
  on the number of cores concurrently employed with high
  efficiency. Moreover, on any given number of cores, the proposed
  parallel approach preserves the $\mathcal{O}(N\log N)$ computing
  cost inherent in the sequential version of the IFGF
  algorithm. Unlike other approaches, the IFGF method does not utilize
  the Fast Fourier Transform (FFT), and it is thus better suited than
  other methods for efficient parallelization in distributed-memory
  computer systems. In particular, the IFGF method relies on a
  ``peer-to-peer'' strategy wherein, at every level, field propagation
  is directly enacted via ``exchanges'' between ``peer'' polynomials
  of low and constant degree, without data accumulation in large-scale
  ``telephone-central'' mathematical constructs such as those in the
  Fast Multipole Method (FMM) or pure FFT-based approaches. A variety
  of numerical results presented in this paper illustrate the
  character of the proposed parallel algorithm, including excellent
  weak and strong parallel scaling properties in all cases
  considered---for problems of up to $4,096$ wavelengths in acoustic
  size, and scaling tests spanning from 1 compute core to all $1,680$
  cores available in the High Performance Computing cluster used.
\end{abstract}
\vspace{0.5 cm}
\maketitle
\section{\label{sec:introduction}Introduction}

This paper presents a parallel implementation of the ``Interpolated
Factored Green Function'' (IFGF) method introduced recently for the
accelerated evaluation of discrete integral operators arising in wave
scattering and other areas~\cite{Bauinger2021}. The proposed
implementation, which is structured as a hybrid MPI-OpenMP computer
program suitable for instantiation in modern high-performance
computing systems (HPC), demonstrates in practice excellent parallel
scaling up to large numbers of cores, without any hard limitations on
the number of cores concurrently employed with high efficiency, while
preserving the linearithmic complexity (namely, $\mathcal{O}(N\log N)$
computing cost) inherent in the sequential IFGF algorithm. The IFGF
method accelerates the evaluation of discrete integral operators by
relying on a certain factorization of the Green function into two
factors, a ``centered factor'' that is incorporated easily as a common
factor in the calculation, and an ``analytic factor'' which enjoys a
property of analyticity up to and including infinity---and which thus
motivates the IFGF strategy, namely, evaluation of a given discrete
integral operator by means of a hierarchical interpolation approach
which relies on use of a large number of small and independent
interpolation procedures. In particular, the IFGF method does not
utilize acceleration elements commonly used by other acceleration
methods~\cite{2007DirectionalFMMLexing, 2012FMMChebyshevMessnerSchanz,
  Cands2009FastButterfly,1996ButterflyMichielssen,ROKHLIN1993,2006FMMRokhlin,2017Boerm,2017BoermH2,2001FFTKunyansky,1996AIMFFTBelszynski,
  1997FFTPhillips,2003FMMLexingKernelIndependent,2003BebendorfRjasanow}
such as the Fast Fourier Transform (FFT), special-function expansions,
high-dimensional linear-algebra factorizations, translation operators,
equivalent sources, or parabolic scaling; detailed discussions in this
regard can be found in~\cite{Bauinger2021}, below in this section, and
indeed, at many points throughout this paper. Roughly speaking, the
IFGF method relies on a ``peer-to-peer'' strategy wherein, at every
level, field propagation is directly enacted via ``exchanges'' between
``peer'' polynomials of low and constant degree, without data
accumulation in large-scale ``telephone-central'' mathematical
constructs which require a ``downward pass'' through the box octree
inherent in the Fast Multipole Method (FMM)~\cite{2006FMMRokhlin}, or
the surface evaluation of equivalent sources in direct FFT-based
methods~\cite{2001FFTKunyansky}. A variety of numerical results
presented in this paper illustrate the character of the proposed
parallel method, including excellent weak and strong parallel scaling
properties in all cases considered---for problems of up to $4,096$
wavelengths in electrical size, and scaling tests spanning from 1
compute core to all $1,680$ cores available in the HPC cluster used.
 
The parallelization of accelerated Green function methods has been the
subject of a significant literature, which is mostly devoted to
tackling a particular difficulty, namely, the ``parallelization
bottleneck''---which manifests itself under various related
guises~\cite{2014ParallelDirectionalFMMLexing,
  ParallelFMMChandramowlishwaran2010, 2021LargeMLFMA, 2019LargeMLFMA,
  2009Ergul, 2007Volakis, 2009LargeSingleLvlFMM, 2014Fanghzou1,
  2014Fangzhou2}, and which almost invariably concerns uses of the
hard-to-parallelize~\cite{Nikolic2014} FFT
algorithm. (Reference~\cite[Sec. 7]{Gumerov2004}, for example,
mentions two alternatives to the use of FFTs in the context of the
FMM, which, however, it discards as less efficient than an FFT-based
procedure.)  In the case of the multilevel Fast Multipole Method
(FMM), the parallelization bottleneck arises in the evaluation of
translation operators associated with the upper part of the octree
structure, which leads to low parallel
efficiency~\cite{2007DirectionalFMMLexing,
  2014ParallelDirectionalFMMLexing, 2007Volakis}. In the
``directional'' FMM~\cite{2014ParallelDirectionalFMMLexing} the low
efficiency in the upper octree is alleviated as a result of the
parabolic scaling utilized; however, the parallelization strategy does
suffer from hard limitations in the number of parallel tasks that, in
the cases considered in that reference, lead to a ``leveling off'' of
the parallel scaling at 256 or 512 cores~\cite[Secs. 3.6,
4.2]{2014ParallelDirectionalFMMLexing}, depending on the geometry
under consideration.
Reference~\cite{ParallelFMMChandramowlishwaran2010} identifies the
part of the FMM relying on FFTs as a parallelization bottleneck which
arises from FFT-related ``lowest arithmetic intensity'' and
``bandwidth contention''.  In \cite{2021LargeMLFMA, 2019LargeMLFMA},
in turn, a hybrid octree storage strategy is used, which stores a
complete set of tree nodes for a certain number of ``full'' levels in
each process, and which reduces the communication in the upper octree
levels. Those articles demonstrate the treatment of problems
containing very large numbers of discretization points on up to
$2,560$ processes, but they restrict their illustration of the
algorithm's parallel efficiency to a limited strong scaling test from
1 process (sequential) to 64 processes. In contrast to this hybrid
octree-storage strategy, reference~\cite{2009Ergul} simultaneously
partitions boxes (clusters) and field values representing the
radiating and incoming fields of each box. This approach leads to
increased efficiency compared to a parallelization purely based on the
boxes (clusters), but the communication in the translation step still
poses a bottleneck, resulting in as little as 30\% parallel efficiency
from one (sequential) core to 128 cores.

Reference~\cite{2019Keyes}, in turn, presents scaling results for the
parallel BEMFMM implementation of the FMM algorithm, for wave
scattering problems on up to 196,608 cores on 6,144 compute nodes, and
for problems with up to 2.3 billion degrees of freedom (DOF) and
approximately 1,389 wavelengths in size (or, in the nomenclature of
Table~2 in~\cite{2019Keyes}, a sphere two-meters in diameter
illuminated at the frequency of $f= 238.086$ KHz, under the assumption
of a 343m/s speed of sound).  Like the implementations mentioned
above, the results in~\cite{2019Keyes} indicate a deterioration of the
strong-scaling for growing numbers of cores, as manifested by a
flattening of the strong-scaling speedup curves presented as the
numbers of cores increase. The weak scaling curve presented
in~\cite{2019Keyes} indicate a high weak-scaling efficiency, however,
with up to 95$\%$ efficiency for weak scaling between 32 and 131,072
cores. Comparison with BEMFMM and IFGF weak scaling results presents
some challenges on a number of counts. On one hand the
contribution~\cite{2019Keyes} does not mention a crucial element in
judging parallelization quality, namely, memory usage: even though
memory duplication may be relied upon in a parallel algorithm to
maximize parallel efficiency, no indications are provided in that
paper about the amount of memory used in any of the runs
presented. Further, under closer examination, the computing times
indicated in these curves appeared to be high, and we thus decided to
perform a direct comparison of the performance of our IFGF
implementation with the BEMFMM implementation on the basis of the
freely available BEMFMM open-source download provided by the
authors. By necessity, our tests were limited to a test example
consisting of a sphere containing approximately 360,000 DOF, which is
the largest test case provided with the BEMFMM test code, and we
selected a sphere of acoustic diameter of $16\lambda$ for this
experiment. We run both algorithms in the 30 available nodes in our
cluster, Wavefield, each one of which contains 56 computing cores. Our
observations are as follows. The BEMFMM run for the test case
considered required 20 secs. in a single node, and 5 secs. in all 30
nodes, with a speedup factor of 4 going from 1 to thirty nodes. The
IFGF run, in turn, required 1.6 secs. in a single node, and 0.122 secs
in the thirty node cluster, with a speedup factor of 13 going from 1
to 30 nodes. Thus, the IFGF runs in one and 30 nodes were faster than
the BEMFMM runs by factors of 12.5 and 40 in the 1-node and 30-node
runs, respectively, with an IFGF speedup over three times higher than
that provided by BEMFMM going from 1 to 30 nodes. As an additional
point of contact with reference~\cite{2019Keyes}, it is worth
mentioning that, in our 1,680 core cluster, and on the basis of
approximately 4 TB of memory, a sphere 1,389 $\lambda$ in diameter
(reported as $f=238.086$KHz at a speed of sound of $343m/s$
in~\cite{2019Keyes}; cf. Section~\ref{sec:largespheretests} for
details) with 2.12 billion DOF was run in a computing time of 2,380
seconds (Table~\ref{table:largesphere} in
Section~\ref{sec:largespheretests}), which, with a 0.5$\%$ near-field
error (which may be compared to the only error indicator reported in
\cite[Table 2]{2019Keyes} for this test case, which amounts to 20$\%$,
as well as the 3$\%$ near-field solution error reported in the same
table of that paper for significantly smaller problems), is  a
factor of approximately 46 times longer than the time reported
in~\cite{2019Keyes}, for the same number of DOFs and sphere size, on a
computer 78 times larger (containing 131,072 cores) and on the basis
of an unspecified amount of memory. Additional test cases for large
sphere problems are presented in Section~\ref{sec:largespheretests}.


Following a different approach, to avoid the communication bottleneck
in the upper multilevel FMM octree entirely,
references~\cite{2007Volakis, 2009LargeSingleLvlFMM} utilize a
single-level Fast Multipole strategy. While this method significantly
simplifies the algorithm and minimizes the required communication in a
parallel setting, it does give rise to a sub-optimal asymptotic
computational cost (e.g. $\mathcal{O}(N^{3/2})$
in~\cite{2009LargeSingleLvlFMM} or, exploiting the FFT,
$\mathcal{O}(N^{4/3} \log^{2/3} N)$ in~\cite{2007Volakis}), and, while
resulting in good parallel scaling up to 512 processes in the
$\mathcal{O}(N^{3/2})$ algorithm~\cite{2009LargeSingleLvlFMM}, as in
the case of \cite{2014ParallelDirectionalFMMLexing}, the parallel
efficiency does level off beyond 512 processes.  Direct FFT methods,
in turn, present alternatives to the various FMM strategies,
including, for example, the Adaptive Integral
Method~\cite{1996AIMFFTBelszynski} (AIM) and the sparse-FFT
method~\cite{2001FFTKunyansky}. Like the single-level FMM algorithms,
these FFT methods exhibit sub-optimal algorithmic complexity (of
orders $\mathcal{O}(N^{3/2})$ and $\mathcal{O}(N^{4/3})$,
respectively, and, owing to their strong reliance on FFTs, they also
suffer from reduced parallel efficiency, as shown and discussed for
the AIM in e.g. ~\cite{2014Fanghzou1, 2014Fangzhou2}. (A parallel
version of the algorithm \cite{2001FFTKunyansky}, which has been
developed by the authors, has not been published, but we report here
that, as may have been expected, the overall parallel efficiency of
the method suffers from the typical FFT-related degradation.)

Finally we mention parallel methods proposed for
non-singular~\cite{2014ButterflyLexing} and
low-frequency~\cite{2012Winkel,2016Biros} problems which, albeit
important and interesting, do not incur some the main challenges
associated with the singular and high-frequency kernels considered in
this paper. Thus, although not applicable to singular Green function
kernels such as the ones considered here, the Butterfly
Method~\cite{2014ButterflyLexing} does provide an acceleration
technique for Fourier integral operators which, based on
linear-algebra constructs instead of the hierarchical interpolation
underlying the IFGF approach, incorporates a parallelization strategy
that is somewhat reminiscent of the proposed IFGF parallelization
approach. The Blue Gene/Q implementation~\cite{2014ButterflyLexing} of
the Butterfly parallel algorithm demonstrates excellent results in
terms of parallel scaling to a large number of cores. The parallel FMM
method presented in~\cite{2016Biros}, which is restricted to box
geometries and to the Laplace and low-frequency Helmholtz problems,
shows impressive scaling up to 299,008 cores on 18,688
nodes. Similarly, the parallel Barnes-Hut tree code~\cite{2012Winkel}
for the low-frequency singular problem provides excellent scaling up
to 294,912 cores with up to 2,048,000,000 particles.

The parallel IFGF strategy introduced in this contribution is based on
adequate partitioning of the interpolations performed on each level of
the underlying octree structure, which facilitates the spatial
decomposition of the surface discretization points. As shown
in~\cite{Bauinger2021}, the number of interpolations performed on each
level is large and approximately constant (as a function of the octree
level). The decomposition and distribution of the interpolation data
is based on a total order in the set of spherical cone segments
representing the interpolation domains, which is an extension of a
domain decomposition based on a Morton curve to the box-cone data
structure inherent in the IFGF approach. The usage of space-filling
curves for the representation of octree structures underlying the
various acceleration methods is not a novel concept~\cite{2019Keyes,
  2013Warren, 2016Biros}. However, the extension of space-filling
curves to the present box-cone structure of the IFGF method to achieve
the desired efficiency has not been reported before. In view of its
strong reliance on the IFGF's box-cone structure, the proposed
parallelization strategy is therefore not applicable to other
acceleration methods such as the FMM. The present parallel IFGF
implementation on a 30-node (1,680-core) HPC cluster with Infiniband
interconnect, delivers perfect $\mathcal{O}(N\log N)$ performance on
all 1,680 cores. And, demonstrating high (albeit imperfect) strong and
weak parallel efficiencies, unlike other methods, it does not suffer
from scaling limitations, under either weak scaling or strong scaling
tests, as the number of processing cores grow---conceivably, as argued
in Sections~\ref{sec:mpiparallelization} and~\ref{sec:results}, up to
arbitrarily large numbers of cores. Note that the presented
implementation of the proposed parallelization approach deliberately
does not utilize any optimizations tailored to the hardware used for
the performance tests. It is therefore reasonable to expect a similar
performance to that presented in this paper, without hardware-specific
optimizations, on a wide variety of modern cluster systems.

This paper includes five sections in its main body
(Sections~\ref{sec:introduction} through~\ref{sec:conclusio}) and two
Supplementary Materials sections (Sections~\ref{sec:hpcnomenclature}
and~\ref{sec:supp_num_res}). Section~\ref{sec:ifgfmethod} briefly
summarizes the description~\cite{Bauinger2021} of the IFGF method, and
it introduces the required notations and
nomenclature. Section~\ref{sec:parallelization} then introduces the
proposed parallelization strategy, including a description of computer
systems considered (in the first paragraphs of
Section~\ref{sec:parallelization}, see also
Section~\ref{sec:hpcnomenclature}), the associated OpenMP and MPI
parallelization approaches (Sections~\ref{sec:ompparallelization}
and~\ref{sec:mpiparallelization}, respectively), and a discussion of
the IFGF's linearithmic communication costs
(Section~\ref{sec:complexity})---which results in preservation of the
overall IFGF linearithmic parallel computing cost. A variety of
numerical results are presented in Section~\ref{sec:results} and, with
additional detail, in the supplementary materials
Section~\ref{sec:supp_num_res}. A few concluding comments, finally,
are presented in Section~\ref{sec:conclusio}.

\section{\label{sec:ifgfmethod}Review of the IFGF Method}
As discussed above, the IFGF method provides an accelerated algorithm,
requiring $\mathcal{O}(N \log{N})$ operations, for the numerical
evaluation of discrete integral operators of the form
\begin{equation} \label{eq:field1}
    I(x_\ell) \coloneqq \sum \limits_{\substack{m = 1 \\ m \neq \ell}}^N a_m G(x_\ell, x_m),\quad \ell = 1, \ldots, N,
\end{equation}
for given points $x_\ell$ on a surface $\Gamma \subset \mathbb{R}^3$, and for given complex coefficients $a_m \in \mathbb{C}$, where the function $G(x, y)$, defined for $x, y \in \mathbb{R}^3$, denotes a Green function for some partial differential equation, such as the acoustic Green function
\begin{equation*}
    G(x, y) = \frac{e^{\iota \kappa |x-y|}}{4 \pi |x-y|}
\end{equation*}
associated with the Helmholtz equation ($\iota$ denotes the imaginary unit and $\kappa$ the wavenumber) as well as those associated with the Laplace, Stokes, and elasticity equations, among others. In what follows we denote by $\Gamma_N := \{x_1, \ldots, x_N\} \subset \Gamma$ the set of surface discretization points.
				
\begin{figure}
   \centering \subcaptionbox{\centering Two-dimensional sketch of a
     three level ($D=3$) IFGF domain decomposition with neighbors (in
     white) and cousin boxes (in gray) for the particular box
     $B_{(2, 1)}^3$. A surrogate scatterer is sketched in
     blue. \label{fig:domaindecomposition}}[0.45\linewidth]{\includegraphics[width=0.4\textwidth]{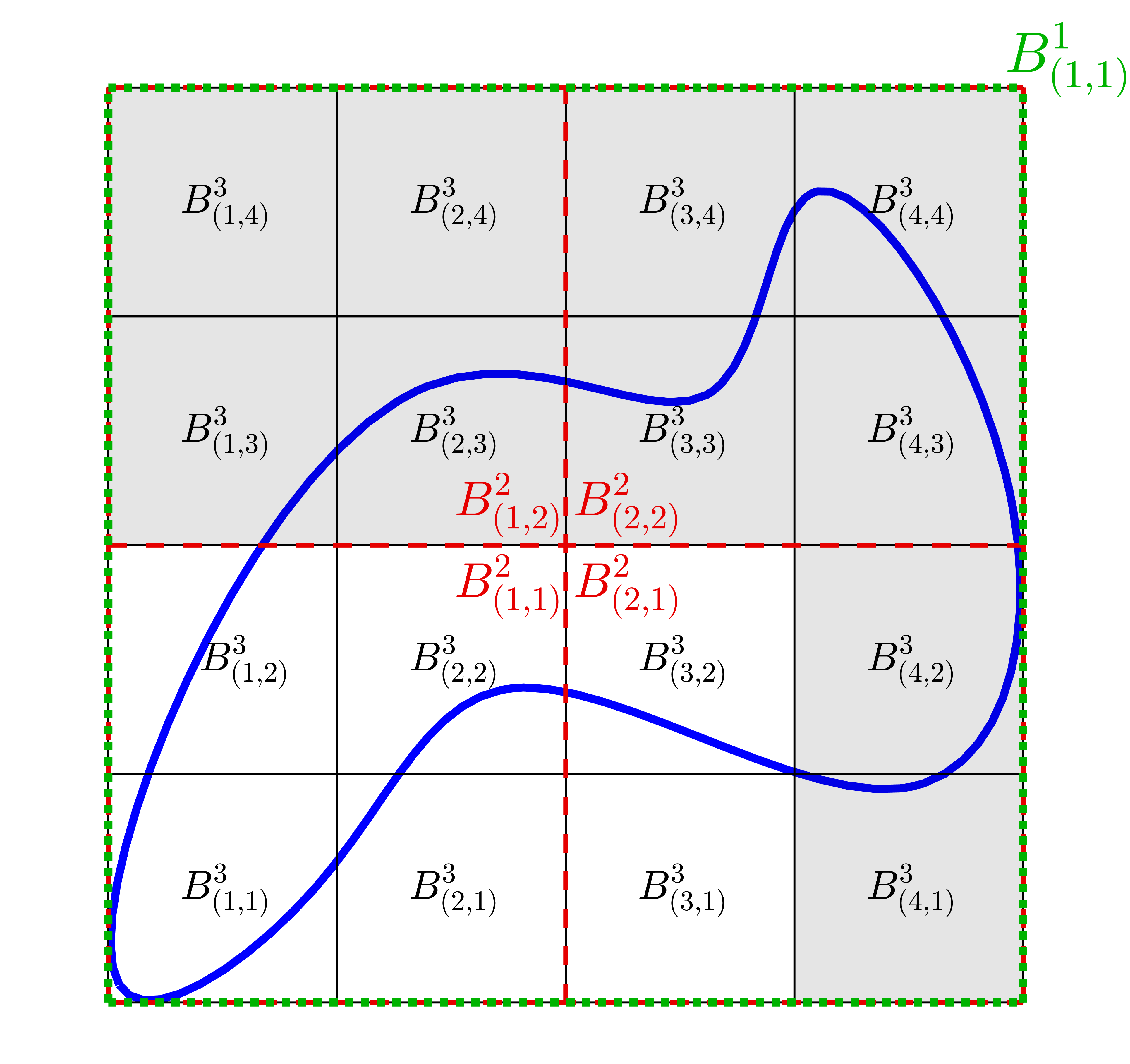}}
   \setbox9=\hbox{\includegraphics[width=.4\linewidth]{multilevel_domain_partition_sketch_2d_cut_2.png}}
   \subcaptionbox{\centering Two-dimensional illustration of cone
     segments: in red, cone-segments co-centered with the box
     $B_\mathbf{k}^d$, and, in black, cone-segments co-centered with
     the parent box
     $\mathcal{P}B_\mathbf{k}^d$.\label{fig:conesegments}}[0.45\linewidth]{\raisebox{\dimexpr\ht9/2-\height/2}{\includegraphics[width=0.4\textwidth]{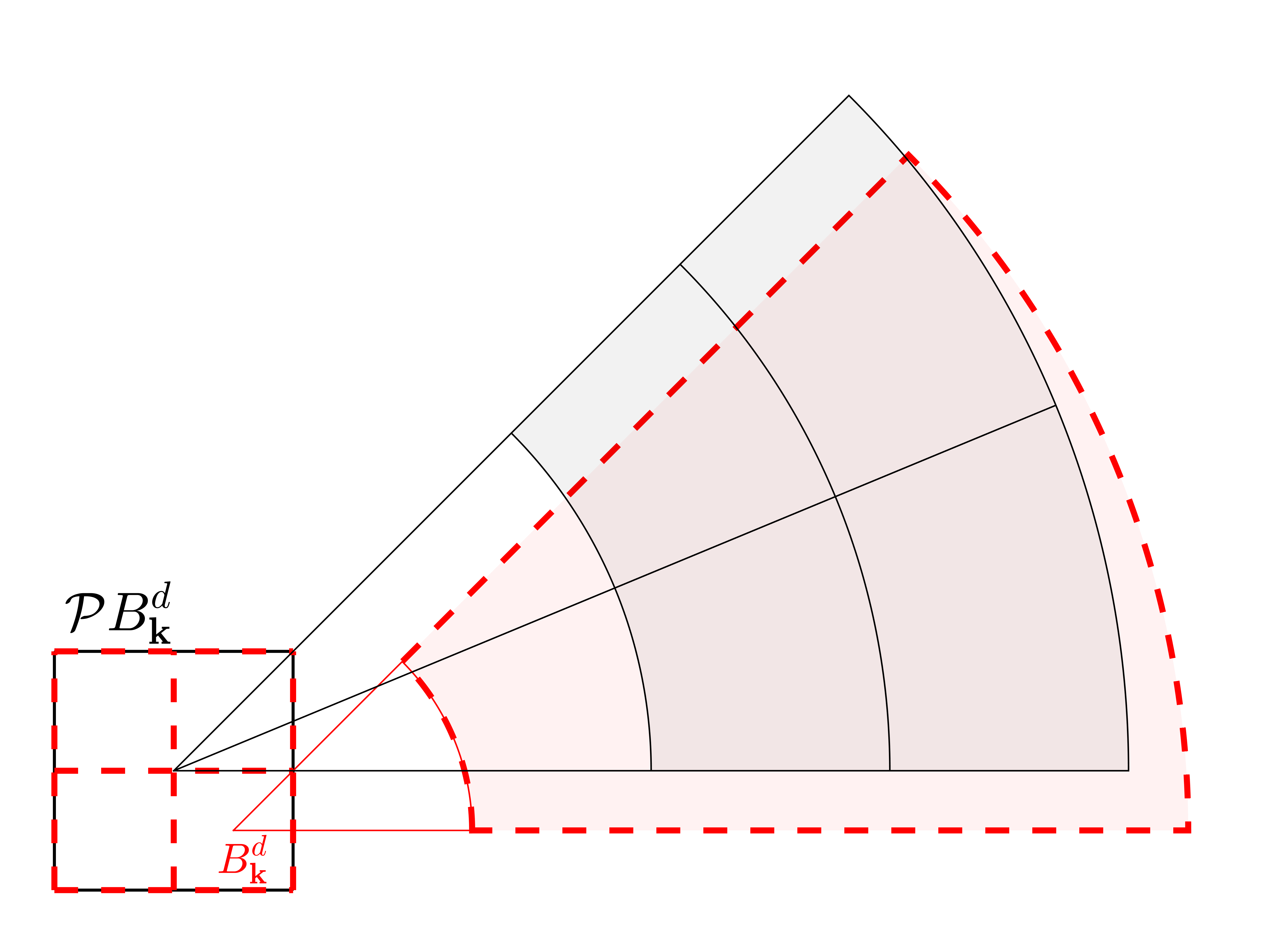}}}
    \caption{Illustration of the box and cone hierarchical structures
      used in the IFGF method.}
  \end{figure}

For a given $D \in \mathbb{N}$, the IFGF method is based upon use of a
$D$-level octree hierarchical decomposition of a cube $B_{1, 1, 1}^1$
containing the discrete surface $\Gamma_N$, where each level is
determined by a set of axis-aligned boxes
$B_\mathbf{k}^d \subset \mathbb{R}^3$ (defined as the Cartesian
product of three one-dimensional half-open intervals of the form
$[a,a+H_d)$ for some $a\in\mathbb{R}$), where
$\mathbf{k} \in \mathbb{N}^3$ denotes a multi-index describing the
three-dimensional position of the box in the resulting Cartesian grid
of boxes, and where $d$ ($1 \leq d \leq D$) denotes the level in the
octree. The octree structure of boxes is defined iteratively starting
from a single box $B_{1, 1, 1}^1 \supset \Gamma_N$ of side
$H_1 \in \mathbb{R}$, $H_1 > 0$, on level $d = 1$. (Note that there is
no undue expense incurred for, say, an elongated surface $\Gamma$, for
which a cubic box could be mostly empty---since, as indicated in what
follows, only certain ``relevant'' child boxes in the box octree are
used by the algorithm.) The boxes on consecutive levels
$d = 2, \ldots, D$ are defined by means of a partition of each of the
level $(d-1)$ boxes into eight equi-sized and disjoint boxes of side
$H_d = H_{d-1}/2$ resulting in the level $d$ boxes $B_\mathbf{k}^d$
($\mathbf{k} \in \{1, \ldots, 2^{d-1}\}^3 =: I^d_B$). The
two-dimensional equivalent of the resulting hierarchical octree
structure for an illustrative three-level configuration ($D = 3$) is
depicted in Figure \ref{fig:domaindecomposition}.
Clearly, each box $B_{\mathbf{k}}^d$ on level $d$ ($2 \leq d \leq D$) admits a unique $(d-1)$-level parent box $\mathcal{P}B_{\mathbf{k}}^d$ containing $B_{\mathbf{k}}^d$.

To achieve the desired acceleration, the IFGF method only considers interactions between boxes in a certain set $\mathcal{R}_B$ of \textit{relevant boxes}, which are defined as the boxes in the octree structure that contain at least one surface discretization point:
\begin{equation*}
\begin{split}
    \mathcal{R}_B^d &:= \{ B_{\mathbf{k}}^d \ : \  \Gamma_N \cap B_{\mathbf{k}}^d \neq \emptyset, \mathbf{k} \in I^d_B \}, \quad  d=1, \ldots, D,  \\
    \mathcal{R}_B &:= \bigcup \limits_{d = 1, \ldots, D} \mathcal{R}_B^d. 
\end{split}
\end{equation*}
Furthermore, for a given box $B_\mathbf{k}^d$ on any level $d$, the method relies on a number of additional concepts, such as the \textit{set of neighboring boxes} $\mathcal{N} B_\mathbf{k}^d$ (namely the set of level-$d$ boxes whose closure have a non-empty intersection with the closure $\overline{B_\mathbf{k}^d}$ of $B_\mathbf{k}^d$)  and the \textit{set of cousin boxes} $\mathcal{M} B_\mathbf{k}^d$ (non-neighboring boxes that are children of the parents neighbors), as well as related concepts such as the \textit{set of neighboring points} $\mathcal{U} B_\mathbf{k}^d$ and the \textit{set of cousin points} $\mathcal{V} B_\mathbf{k}^d$ (which denote the set of surface discretization points within the neighboring boxes and the cousin boxes, respectively):
\begin{equation} \label{eq:defneighboursandcousins}
\begin{split}
    \mathcal{N} B_\mathbf{k}^d &:= \{ B_\mathbf{j}^d \in \mathcal{R}_B : ||\mathbf{j} - \mathbf{k}||_\infty \leq 1 \}, \\
    \mathcal{M} B_\mathbf{k}^d &:= \{ B_\mathbf{j}^d \in \mathcal{R}_B : B_\mathbf{j}^d \notin \mathcal{N} B_\mathbf{k}^d \wedge \mathcal{P}B_{\mathbf{j}}^d \in \mathcal{N} \mathcal{P} B_\mathbf{k}^d \}, \\
    \mathcal{U} B_\mathbf{k}^d &:= \left( \bigcup \limits_{B \in \mathcal{N} B_\mathbf{k}^d} B \right) \cap \Gamma_N,  \\
    \mathcal{V} B_\mathbf{k}^d &:= \left( \bigcup \limits_{B \in \mathcal{M} B_\mathbf{k}^d} B \right) \cap \Gamma_N.
\end{split}
\end{equation}
Figure \ref{fig:domaindecomposition} displays the neighbors and cousins of the box $B_{(2, 1)}^3$ in white and gray colors, respectively.

The IFGF algorithm accelerates the evaluation of the operator
\eqref{eq:field1} by exploiting a fast method for computations of
pairwise interactions between cousin boxes on every level $d$, for
$d = D, \ldots, 3$. (Note that the algorithm ends at level
$d= 3$---since at $d=3$, each box is a cousin or a neighbor of all the
other boxes, and thus, all remaining surface evaluations are completed
at this stage.) The evaluation of these interactions is enacted by
means of a simple piecewise interpolation method based on a certain
factored form of the Green function in a set of box-centered spherical
coordinate systems, with one such spherical-coordinate system centered
at each one of the relevant boxes. The use of angular and radial
interpolation methods gives rise to so-called \textit{cone segments}
$C_{\mathbf{k}; \bf{\gamma}}^d$: one for each box $B_\mathbf{k}^d$ and
for each multi-index ${\bf \gamma} \in I_C^d \subset \mathbb{N}^3$
characterizing a conical interpolation domain. A set
$\mathcal{X} C_{\mathbf{k}; \bf{\gamma}}^d$ of $P \in \mathbb{N}$
\textit{interpolation points} is used within each cone segment
$C_{\mathbf{k}; \bf{\gamma}}^d$, where $P$ is an arbitrary but fixed
number throughout the algorithm. In what follows, cone segments
$C_{\mathbf{k}; \bf{\gamma}}^d$ are called \textit{co-centered} with a
box $B_\mathbf{k}^d$ if and only if the origin of the spherical
coordinate system defining the cone segment coincides with the center
of the box. Two cone segments are called \textit{co-centered} if they
are co-centered with the same box. Note that the sub- and super
indices $\mathbf{k}$ and $d$ in the cone-segment notation
$C_{\mathbf{k}; \bf{\gamma}}^d$ coincide with the corresponding
indices of the co-centered box.

To achieve competitive computation times, a certain factorization of the Green function $G(x, x') = G(x, x_\mathbf{k}^d) g_\mathbf{k}^d(x, x')$ into a \textit{centered} factor $G(x, x_\mathbf{k}^d)$ (centered at the
\textit{box-center} $x_\mathbf{k}^d$ of the box $B_\mathbf{k}^d$) and
an \textit{analytic} factor $g_\mathbf{k}^d(x, x')$, is used. The
field $I(x)$ in \eqref{eq:field1} can be expressed, for each level
$d$, as the sum, over all multi-indices $\mathbf{k} \in I_B^d$, of
fields $I_\mathbf{k}^d(x)$ equal to the sum of $G(x,x')$ for all
surface discretization points $x'$ within the box $B_\mathbf{k}^d$,
i.e., for all $x' \in B_\mathbf{k}^{d} \cap \Gamma_N$. Using the
aforementioned factorization centered at $x_\mathbf{k}^d$ yields
\begin{equation}\label{eq:fieldperbox}
    I_\mathbf{k}^d(x) = \sum \limits_{x' \in B_\mathbf{k}^{d} \cap \Gamma_N} a(x') G(x, x') = G(x, x_\mathbf{k}^d) F_\mathbf{k}^d (x),  \qquad F_\mathbf{k}^d(x) \coloneqq \sum \limits_{x' \in B_\mathbf{k}^{d} \cap \Gamma_N} a(x') g_\mathbf{k}^d(x, x'),
\end{equation}
where $a(x')$ denotes the coefficient $a_m$ in \eqref{eq:field1} that
corresponds to the point $x' \in \Gamma_N$. The IFGF interpolation
procedure is then used to evaluate $F_\mathbf{k}^d$. The generation of
the $P$ coefficients of each one of the degree-$P$ polynomial
interpolants $I_P C_{\mathbf{k}; \bf{\gamma}}^d$, corresponding to
interpolation of the field $F_\mathbf{k}^d$ (cf. \eqref{eq:fieldperbox})
over the cone segment $C_{\mathbf{k}; \bf{\gamma}}^d$, is achieved on
the basis of the field values $F_\mathbf{k}^d(\mathcal{X} C_{\mathbf{k}; \bf{\gamma}}^d) := \{
F_\mathbf{k}^d(x) : x \in \mathcal{X} C_{\mathbf{k}; \bf{\gamma}}^d \}$.

In \cite{Bauinger2021} it is shown that the analytic factor is
analytic everywhere in
$\mathbb{R}^3\setminus \mathcal{N} B_\mathbf{k}^d$ and up to and
including infinity, and it can therefore be interpolated accurately
throughout that region on the basis of a small (finite!) number of
interpolation points. (It is is easy to check that the same is true
for most of the relevant kernels arising in applications.) Since
$F_\mathbf{k}^d$ equals a linear combination of finitely many
analytic-factor functions, it is clear that this function shares the
same analytic properties, and it can therefore be interpolated with
equal quality and efficiency.  The cone segments
$C_{\mathbf{k}; \bf{\gamma}}^d$ are defined by means of an iterative
procedure similar to the one used in the definition of the boxes
$B_\mathbf{k}^d$, but in reversed order, starting from $d=D$ and moving
upwards the tree to $d=3$.  The set of cone segments that is to be
used at a given level $d$ depends strongly on the character of the
surface $\Gamma_N$, the wavenumber $\kappa$ and, possibly, the Green
function $G$. A two-dimensional sketch of some illustrative
box-centered cone segments for a given box $B_\mathbf{k}^d$ and its
parent $\mathcal{P}B_\mathbf{k}^d$ is provided in Figure
\ref{fig:conesegments}.

\begin{figure}
    \centering
    \includegraphics[width=.9\linewidth]{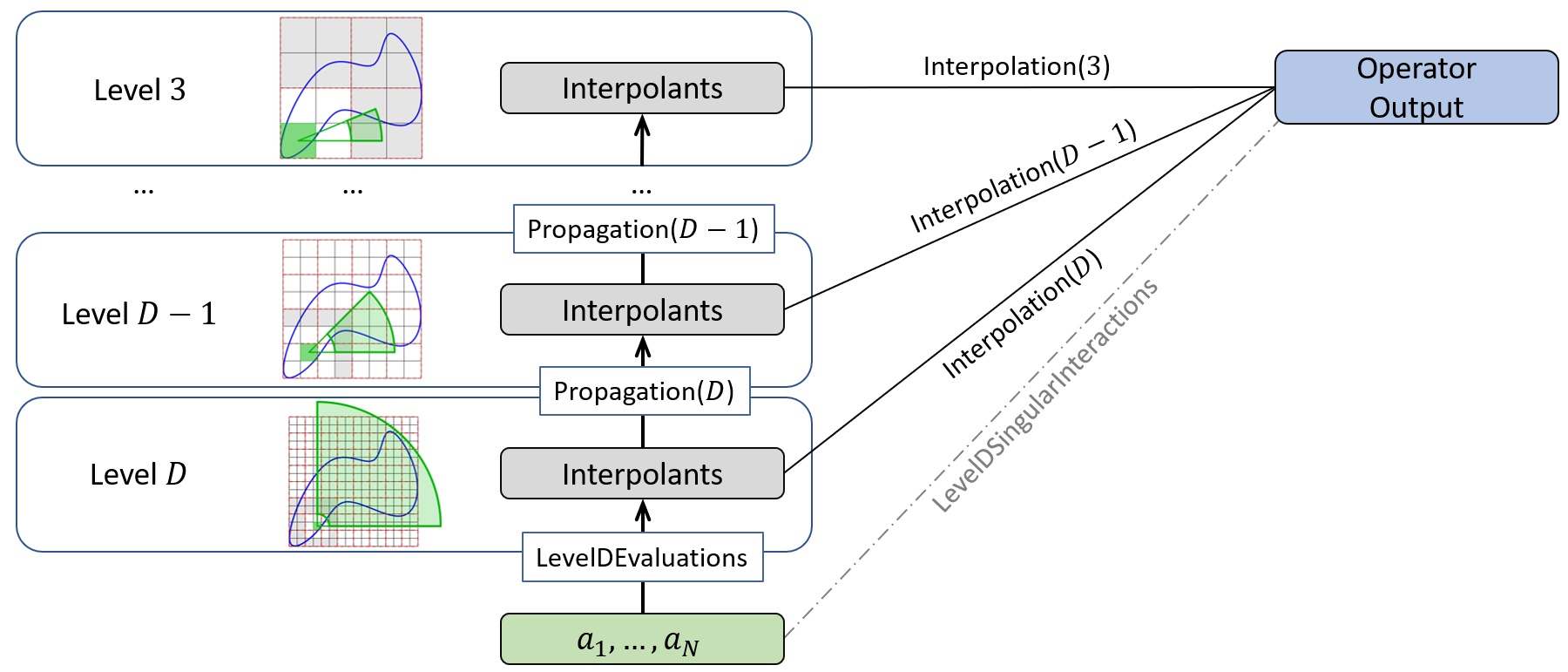}
    \caption{Visual representation of the IFGF algorithm, outlined in
      Algorithm~\ref{alg:ifgf}, and also expressed in
      Algorithm~\ref{alg:ifgf_short} in terms of three fundamental
      functions called \textit{LevelDEvaluations},
      \textit{Propagation} and \textit{Interpolation}. Starting from
      the given coefficients $a_1, \ldots, a_N$ in
      equation~\eqref{eq:field1}, the \textit{LevelDEvaluations}
      function generates the first set of interpolants on level
      $D$. The \textit{Interpolation} function interpolates to the
      surface discretization points and the \textit{Propagation}
      function facilitates the upwards traversal of the octree
      structure. Although they are not part of the IFGF algorithm, the
      interactions between level-$D$ neighbor boxes are represented
      here by the \textit{LevelDSingularInteractions} function. Note
      that, unlike other acceleration methods such as the FMM,
      contributions to the operator output are made at every level,
      and without a requirement of a downward pass over the octree.}
    \label{fig:sketchifgfflow}
  \end{figure}

In order to evaluate the discrete operator~\eqref{eq:field1} in
$\mathcal{O}(N \log N)$ operations, the IFGF algorithm uses iterated
interpolation, as illustrated in Figure~\ref{fig:sketchifgfflow}, to
evaluate the analytic factor at the interpolation points of
consecutive levels---thus avoiding the cost of re-evaluating the field
$I_\mathbf{k}^d(x)$ on levels $(d-1),(d-2),\dots, 3$, and using
instead the interpolation data on level $d$ to generate the necessary
interpolation data at on the consecutive level $(d-1)$.  It is
important to note that, in order to further increase the efficiency
and achieve the desired $\mathcal{O}(N\log N)$ complexity, in analogy
to the approach used for boxes, the IFGF method only utilizes the set
of \textit{relevant cone segments} $\mathcal{R}_C B_\mathbf{k}^d$ for
each box $B_\mathbf{k}^d$, namely, the cone segments that are actually
needed for interpolation back to cousin surface discretization points
or to relevant cone segments on the parent level. In other words, the
relevant cone segments $\mathcal{R}_C B_\mathbf{k}^d$ are defined by
\begin{equation}\label{eq:relevantconesegments}
  \begin{split}
    \mathcal{R}_C B_\mathbf{k}^d &:= \emptyset \quad \text{for} \quad d = 1,2,\\
    \mathcal{R}_C B_\mathbf{k}^d &:= \left\{ C_{\mathbf{k}; \mathbf{\gamma}}^d \, : \, \gamma \in I_C^d\, , \, C_{\mathbf{k}; \mathbf{\gamma}}^d \cap \mathcal{V} B_\mathbf{k}^{d} \neq \emptyset \text{ or } C_{\mathbf{k}; \mathbf{\gamma}}^d \cap \left( \bigcup \limits_{C \in \mathcal{R}_C \mathcal{P} B_\mathbf{k}^{d}}  C \right) \neq \emptyset \right \} \quad \text{for} \quad d \geq 3.
\end{split}
\end{equation}

The serial IFGF algorithm, introduced in \cite{Bauinger2021}, is
  summarized in Figure~\ref{fig:sketchifgfflow} and described in what
  follows. Starting from the given coefficients in
  equation~\eqref{eq:field1} at the bottom of
  Figure~\ref{fig:sketchifgfflow}, the IFGF algorithm first performs
  ``LevelDSingularInteractions'' (which, while required for the full
  evaluation of~\eqref{eq:field1}, are not, strictly speaking, a part
  of the IFGF strategy, and would, in the context of a scattering
  solver, be substituted by an appropriate local integration scheme;
  see e.g.~\cite{jimenez2021ifgf}). This
  ``LevelDSingularInteractions'' stage evaluates the field
  $I^D_\mathbf{k}$ at all surface discretization points $x$ in all the
  neighbor boxes of each box $B^D_\mathbf{k}$ (i.e. at all
  $x \in \mathcal{U} B^D_\mathbf{k}$ for all relevant boxes
  $B^D_\mathbf{k}$). Next, the algorithm performs
  ``LevelDEvaluations'', that is, it first evaluates the field
  $F^D_\mathbf{k}(x)$ (see \eqref{eq:fieldperbox}) at every
  interpolation point in the relevant cone segments co-centered with
  the box $B^D_\mathbf{k}$, and it subsequently generates the
  necessary level-$D$ interpolants. The IFGF algorithm then proceeds
  through levels $d = D, \ldots, 3$ by performing, on each level $d$,
  1)~Interpolations to cousin surface discretization points in the
  ``Interpolation'' stage, as well as, 2)~Interpolations to level
  $d-1$ interpolation points, and subsequent generation of the
  interpolants on level $d-1$, in the ``Propagation'' stage, just as
  in the ``LevelDEvaluations'' stage, but utilizing the interpolants
  instead of direct field evaluations.

  The corresponding pseudo-code is presented in Algorithm
  \ref{alg:ifgf}. Note that this algorithm does not include
  evaluations of interactions between neighboring boxes on the lowest
  level $D$ (``LevelDSingularInteractions'' in
  Figure~\ref{fig:sketchifgfflow}), which would generally be produced
  by means of a separate algorithm, as mentioned above in this
  section.
 
 \begin{algorithm}
  \small
\begin{algorithmic}[1]
  \State \textbackslash \textbackslash Direct evaluations on the lowest level. 
  \For{$B_\mathbf{k}^D \in \mathcal{R}_B$} \label{algstate:looprelevantboxeslevelD} 
  \For{$C_{\mathbf{k}; \mathbf{\gamma}}^D \in \mathcal{R}_C B_\mathbf{k}^D$} \label{algstate:looprelevantconeslevelD} \Comment{Evaluate $F$ at all relevant interpolation points} 
  \State Evaluate and store $F_\mathbf{k}^D(\mathcal{X} C_{\mathbf{k};\mathbf{\gamma}}^D)$ \label{algstate:levelDEval}
  \State Generate interpolant $I_P C_{\mathbf{k}; \mathbf{\gamma}}^D$
  \EndFor 
  \EndFor
  \State
  \State \textbackslash \textbackslash Interpolation onto surface discretization points and parent
  interpolation points.
  \For{$d = D, \ldots, 3$}\label{algstate:loopd}
  \For{$B_\mathbf{k}^d \in \mathcal{R}_B$} 
  \label{algstate:looprelevantboxes} 
  \For{$x \in \mathcal{V} B_\mathbf{k}^d$} \label{algstate:loopnearestneighbouringsurfacepoints}
  \Comment{Interpolate at cousin surface points} 
  \State Determine $C_{\mathbf{k};\mathbf{\alpha}}^{d}$ s.t. $x \in C_{\mathbf{k};\mathbf{\alpha}}^{d}$ 
  \State Evaluate and add to result $I_P C_{\mathbf{k};\mathbf{\alpha}}^{d}(x) \times G(x, x_\mathbf{k}^d)$ \label{algstate:interpolationtosurfacepoints}
  \EndFor 
  \If {$d > 3$}
  \Comment{Evaluate $F$ on parent interpolation points} 
  \State Determine parent $B_\mathbf{j}^{d-1} = \mathcal{P} B_\mathbf{k}^d$
  \For{$C_{\mathbf{j}; \mathbf{\gamma}}^{d-1} \in \mathcal{R}_C B_\mathbf{j}^{d-1}$} \label{algstate:looprelevantcones} 
  \For{$x \in \mathcal{X} C_{\mathbf{j};\mathbf{\gamma}}^{d-1}$} \label{algstate:loopinterppoints}
  \State Determine $C_{\mathbf{k};\mathbf{\alpha}}^{d}$ s.t. $x \in C_{\mathbf{k};\mathbf{\alpha}}^{d}$ 
  \State Evaluate and add $I_P C_{\mathbf{k};\mathbf{\alpha}}^{d}(x) \times G(x, x_\mathbf{k}^d)/ G(x, x_\mathbf{j}^{d-1})$
  \EndFor 
  \EndFor 
  \EndIf 
  \EndFor
  \Comment{Generate interpolants on parent level} 
  \For{$B_\mathbf{j}^{d-1} \in \mathcal{R}_B$} 
  \For{$C_{\mathbf{j}; \mathbf{\gamma}}^{d-1} \in \mathcal{R}_C B_\mathbf{j}^{d-1}$}
  \State Generate interpolant $I_P C_{\mathbf{j}; \mathbf{\gamma}}^{d-1}$
  \EndFor
  \EndFor
  \EndFor
  \caption{IFGF Method}
\label{alg:ifgf}
\end{algorithmic}
\end{algorithm}

\section{Parallel IFGF Method} \label{sec:parallelization} The IFGF
parallelization scheme proposed in this paper relies on use of a
hybrid MPI-OpenMP approach. As detailed in
Section~\ref{sec:mpiparallelization}, the MPI interface plays two
distinct roles in the proposed scheme: it is used to 1)~Enable
distributed-memory parallelization across compute nodes, and 2)~In the
particular case in which MPI ranks are pinned to NUMA nodes
(non-uniform memory access), to distribute work and handle memory
access across NUMA nodes within each compute node. Additional details
concerning the architecture of the computer used, and, in particular,
NUMA nodes, can be found in Section~\ref{sec:hpcnomenclature}.  The
strategy in point~2) guarantees that memory held by a certain MPI rank
is stored within a single NUMA node and can therefore be accessed
quickly by all cores within the NUMA node. Moreover, in case 2),
access to memory in a different NUMA node within the same compute node
is algorithmically effected through MPI in the same manner as access
to memory in a different compute node.

The OpenMP parallelization, described in Section
\ref{sec:ompparallelization}, is used to further distribute the work
assigned to each MPI rank.  Hence, in the specific hardware
implementation demonstrated in this paper (which is based on use of
compute nodes containing four fourteen-core NUMA nodes) typically four
intra-node MPI ranks are used per compute node, each pinned to a
single NUMA node,  each one of which spawns fourteen OpenMP
threads---which, according to our experiments, leads to the best
performance achievable without the adverse impact (on e.g. code
complexity, memory requirements, or communication) entailed in pure
MPI parallelism within each node.
A general discussion on the performance of hybrid MPI-OpenMP
approaches can be found in~\cite{Bull2010, Drosinos2004,
  Akhmetova2007}.

\subsection{IFGF OpenMP parallelization} \label{sec:ompparallelization} Before introducing
the proposed OpenMP parallelization scheme, we briefly consider a
straightforward OpenMP parallelization strategy which we do not
recommend but which we present for reference. This straightforward and
easily implemented strategy results by simply implementing the
algorithm depicted in Figure~\ref{fig:sketchifgfflow} by assigning,
at each level $d$, the work associated with groups of relevant boxes
to various OpenMP threads (e.g. with a ``\#pragma omp parallel for''
statement), in such a way that each group is handled by a single
thread.  Equi-distribution of relevant boxes onto the OpenMP threads
implies equi-distribution of both the surface discretization points
and the computations performed per thread---but only provided 1)~The
surface discretization points are roughly equi-distributed among the
relevant boxes, and, 2)~There is a sufficient number of relevant boxes
to occupy all OpenMP threads.  The difficulties associated with
point~1) could be negotiated, in view of the law of large numbers
\cite[Section 13]{dekking2005modern}, provided sufficiently many boxes
are used, that is to say, provided point 2)~is satisfied. In other
words, the feasibility of the approach under consideration hinges on
the existence of sufficiently many relevant boxes on each level, as
required by point~2). Unfortunately, however, for any given
discretized surface $\Gamma_N$, point~2) is not satisfied at certain
levels $d$ in the octree structure, unless only a small number of
threads is employed.  Noting that, for any surface $\Gamma_N$, there
are only sixty-four boxes overall on level $d=3$ of the algorithm
(and, in general, even fewer relevant boxes), we see that a definite
limit exists on the parallelism achievable by this approach.  The
method presented in~\cite{2014ParallelDirectionalFMMLexing} uses this
strategy in an MPI context, and it is therefore subject to such a hard
limitation on achievable parallelism (although in a somewhat mitigated
form, owing to the characteristics of that algorithm, as discussed in
Section~\ref{sec:introduction}).  To avoid such limitations we
consider an alternate OpenMP parallelization strategy specifically
enabled by the characteristics of the IFGF algorithm, as described in
what follows.

\begin{algorithm}[ht]
\small
\begin{algorithmic}[1]
  \For{$B_\mathbf{k}^D \in \mathcal{R}_B$}
  \For{$C_{\mathbf{k}; \mathbf{\gamma}}^D \in \mathcal{R}_C B_\mathbf{k}^D$}
  \State Evaluate and store $F_\mathbf{k}^D(\mathcal{X} C_{\mathbf{k};\mathbf{\gamma}}^D)$ 
  \State Generate interpolant $I_P C_{\mathbf{k}; \mathbf{\gamma}}^D$
  \EndFor 
  \EndFor
  \caption{LevelDEvaluations}
\label{alg:ifgf_levelD}
\end{algorithmic}
\end{algorithm}

\begin{algorithm}[ht]
\small
\begin{algorithmic}[1]
  \For{$B_\mathbf{k}^d \in \mathcal{R}_B$}
  \For{$x \in \mathcal{V} B_\mathbf{k}^d$}
  \State Determine $C_{\mathbf{k};\mathbf{\alpha}}^{d}$ s.t. $x \in C_{\mathbf{k};\mathbf{\alpha}}^{d}$ 
  \State Evaluate and add to result $I_P C_{\mathbf{k};\mathbf{\alpha}}^{d}(x) \times G(x, x_\mathbf{k}^d)$ 
  \EndFor 
  \EndFor 
  \caption{Interpolation$(d)$}
\label{alg:ifgf_interpolation}
\end{algorithmic}
\end{algorithm}

\begin{algorithm}[ht]
\small
\begin{algorithmic}[1]
  \For{$B_\mathbf{k}^d \in \mathcal{R}_B$}
  \State Determine parent $B_\mathbf{j}^{d-1} = \mathcal{P} B_\mathbf{k}^d$
  \For{$C_{\mathbf{j}; \mathbf{\gamma}}^{d-1} \in \mathcal{R}_C B_\mathbf{j}^{d-1}$}
  \For{$x \in \mathcal{X} C_{\mathbf{j};\mathbf{\gamma}}^{d-1}$}
  \State Determine $C_{\mathbf{k};\mathbf{\alpha}}^{d}$ s.t. $x \in C_{\mathbf{k};\mathbf{\alpha}}^{d}$ 
  \State Evaluate and add $I_P C_{\mathbf{k};\mathbf{\alpha}}^{d}(x) \times G(x, x_\mathbf{k}^d)/ G(x, x_\mathbf{j}^{d-1})$
  \EndFor 
  \EndFor 
  \EndFor
  \For{$B_\mathbf{j}^{d-1} \in \mathcal{R}_B$} 
  \For{$C_{\mathbf{j}; \mathbf{\gamma}}^{d-1} \in \mathcal{R}_C B_\mathbf{j}^{d-1}$}
  \State Generate interpolant $I_P C_{\mathbf{j}; \mathbf{\gamma}}^{d-1}$
  \EndFor
  \EndFor
  \caption{Propagation$(d)$}
\label{alg:ifgf_propagation}
\end{algorithmic}
\end{algorithm}

\begin{algorithm}[ht]
\small
\begin{algorithmic}[1]
  \State LevelDEvaluations()
  \State
  \For{$d = D, \ldots, 3$}
  \State Interpolation(d)
  \If {$d > 3$}
  \State Propagation(d)
  \EndIf 
  \EndFor
  \caption{IFGF Method}
\label{alg:ifgf_short}
\end{algorithmic}
\end{algorithm}

The proposed strategy proceeds via parallelization of the three
independent programming functions that comprise the IFGF method,
namely the \textit{LevelDEvaluations} function, the
\textit{Interpolation} function and the \textit{Propagation} function,
as mentioned in Section~\ref{sec:ifgfmethod} and illustrated in
Figure~\ref{fig:sketchifgfflow}. The first of these functions, the
\textit{LevelDEvaluations} function, which corresponds to the loop in
line \ref{algstate:looprelevantboxeslevelD} of
Algorithm~\ref{alg:ifgf}, evaluates, for each relevant level-$D$ box,
the field generated by the point sources within the box (given by
\eqref{eq:fieldperbox} with $d=D$) at the interpolation points in all
relevant cone segments co-centered with the box and generates the
required interpolants. The second function, the level-$d$-dependent
\textit{Interpolation} function, which corresponds to line 14 under
the loops in lines 12 and 13 and is represented in
Figure~\ref{fig:sketchifgfflow} by rightward lines connecting various
levels to the ``Operator Output'', performs the necessary
interpolations to cousin surface discretization points on level $d$
($d = 3, \ldots, D$).  The third and final programming function, the
level-$d$-dependent \textit{Propagation} function, which corresponds
to line 20 under the loops in lines 12, 18, and 19 and is represented
in Figure~\ref{fig:sketchifgfflow} by means of upward pointing arrows
targeting the ``Interpolant'' boxes, interpolates, for each relevant
level-$d$ box, to interpolation points in the relevant cone segments
co-centered with the parent box on level $(d-1)$ and generates the
required interpolants. These three functions are outlined in
Algorithms \ref{alg:ifgf_levelD}, \ref{alg:ifgf_interpolation}, and
\ref{alg:ifgf_propagation}, respectively. Using these three functions,
the IFGF algorithm (Algorithm~\ref{alg:ifgf}) may be re-expressed as
Algorithm \ref{alg:ifgf_short}. In what follows, we present our
strategies for efficient parallelization of each one of these
functions separately.

Our approach for an efficient parallelization of the
\textit{LevelDEvaluations} function is based on changing the viewpoint
from iterating through the level-$D$ relevant boxes to iterating
through the set $\mathcal{R}_C^D$ of \textit{all relevant cone
  segments on level $D$}. Since corresponding sets of level-$d$
relevant cone segments for the wider range $3\leq d\leq D$ are
utilized in the parallelization of the \textit{Propagation} function,
we generalize the definition: the set of \textit{all relevant cone
  segments on level $d$} is denoted by $\mathcal{R}_C^d$, that is
\begin{equation} \label{eq:relconesegmentslevelD} \mathcal{R}_C^d :=
  \bigcup \limits_{\mathbf{k} \in I_B^d: B_\mathbf{k}^d \in
    \mathcal{R}_B} \mathcal{R}_C B_\mathbf{k}^d, \quad \text{for } 3
  \leq d \leq D.
\end{equation}
Using \eqref{eq:relconesegmentslevelD}, a parallel version of
Algorithm \ref{alg:ifgf_levelD} is presented in Algorithm
\ref{alg:ifgf_parallellevelD}. The aforementioned change in viewpoint
corresponds to collapsing the two outermost nested loops in Algorithm
\ref{alg:ifgf_levelD}, effectively increasing the number of
independent tasks and, consequently, the achievable parallelism. Note
that, in a C++ implementation, the ``parallel for'' construct in
Algorithm \ref{alg:ifgf_parallellevelD} corresponds to e.g. a ``for''-loop
preceded by the pragma directive ``omp parallel for''.
\begin{algorithm}[tb]
\small
\begin{algorithmic}[1]
  \ParFor{$C_{\mathbf{k}; \mathbf{\gamma}}^D \in \mathcal{R}_C^D$}
  \State Evaluate and store $F_\mathbf{k}^D(\mathcal{X} C_{\mathbf{k};\mathbf{\gamma}}^D)$
  \State Generate interpolant $I_P C_{\mathbf{k}; \mathbf{\gamma}}^D$
  \EndParFor
  \caption{Parallel LevelDEvaluations}
\label{alg:ifgf_parallellevelD}
\end{algorithmic}
\end{algorithm}

The proposed parallelization of the $d$-dependent \textit{Propagation}
function follows a similar idea as the parallel
\textit{LevelDEvaluations} considered above---relying now on iteration
over the relevant $(d-1)$ (parent-level) cone segments, which are
targets of the interpolation, instead of the relevant level-$d$ boxes
emitting the field. This strategy addresses the difficulties arising
from the straightforward approach described at the beginning of
Section~\ref{sec:ompparallelization}, for which the number of
available tasks to be distributed decreases with the level $d$ and
imposes a hard limit on the achievable parallelism. Indeed, in the
context of the oscillatory Green functions over two-dimensional
surfaces $\Gamma \subset \mathbb{R}^3$ considered in this paper, for
example, wherein the number of relevant cone segments on each level is
an approximately constant function of $d$~\cite[Sec. 3.3.3]{Bauinger2021},
the number of independent tasks available for parallelization remains
approximately constant as a function of $d$.  Additionally, the
proposed parallel \textit{Propagation} strategy avoids a significant
``thread-safety'' \cite{reinders2021data, pacheco2011introduction}, predicament, that is ubiquitous in the
straightforward approach, whereby multiple writes to the same target
interpolation point on the parent level take place from different
threads. 
In contrast, the proposed \textit{Propagation} strategy,
is by design thread-safe without any additional considerations, since
it distributes the targets of the interpolation to the available
threads.  Note that the practical implementation of this approach
requires the algorithm to first determine the relevant box
\begin{equation}
  \mathcal{R}_B C_{\mathbf{k}; \bf{\gamma}}^d := B_\mathbf{k}^d \qquad \text{ s.t. }  \quad C_{\mathbf{k}; \bf{\gamma}}^d \in \mathcal{R}_C B_\mathbf{k}^d
\end{equation}
that is co-centered with a given relevant cone segment
$C_{\mathbf{k}; \bf{\gamma}}^d$; then to determine the relevant
level-$(d+1)$ \textit{child boxes}
\begin{equation}
  \mathcal{C}B_\mathbf{k}^d := \left\{ B_\mathbf{j}^{d+1} \in \mathcal{R}_B : \mathbf{j} \in I_B^{d+1}, \mathcal{P} B_\mathbf{j}^{d+1} = B_\mathbf{k}^d \right\}, \end{equation}
of a given relevant box $B_\mathbf{k}^d$ on level $d$; and, finally, to find all the interpolants $I_P C_{\mathbf{k}, \gamma}^d$ on 
the relevant cone segments~\eqref{eq:relevantconesegments} co-centered with the child boxes from which the propagation needs to be enacted.
Using this notation, the resulting \textit{Parallel Propagation} algorithm is presented in Algorithm \ref{alg:ifgf_propagation_parallel}.

\begin{algorithm}[tb]
\small
\begin{algorithmic}[1]
  \ParFor{$C_{\mathbf{j}; \mathbf{\gamma}}^{d-1} \in \mathcal{R}_C^{d-1}$}
  \For{$B_\mathbf{k}^d \in \mathcal{C} (\mathcal{R}_B C_{\mathbf{j}; \mathbf{\gamma}}^{d-1})$}
  \For{$x \in \mathcal{X} C_{\mathbf{j};\mathbf{\gamma}}^{d-1}$}
  \State Determine $C_{\mathbf{k};\mathbf{\alpha}}^{d}$ s.t. $x \in C_{\mathbf{k};\mathbf{\alpha}}^{d}$ 
  \State Evaluate and add $I_P C_{\mathbf{k};\mathbf{\alpha}}^{d}(x) \times G(x, x_\mathbf{k}^d)/ G(x, x_\mathbf{j}^{d-1})$
  \EndFor 
  \EndFor 
  \State Generate interpolant $I_P C_{\mathbf{j}; \mathbf{\gamma}}^{d-1}$
  \EndParFor
  \caption{Parallel Propagation$(d)$}
\label{alg:ifgf_propagation_parallel}
\end{algorithmic}
\end{algorithm}

The proposed parallelization strategy for the third and final IFGF
programming function, namely, the \textit{Interpolation} function,
relies once again on the strategy used for the
\textit{LevelDEvaluations} and \textit{Propagation} functions---which,
in the present case, leads to changing the viewpoint from iterating
through the relevant boxes to iterating through the surface
discretization points that are the target of the interpolation
procedure. This approach avoids both, the difficulties mentioned at
the beginning of Section~\ref{sec:ompparallelization} (concerning the
existence of a small number of relevant boxes in the upper levels of
the octree structure), as well as thread-safety difficulties similar
to those discussed above in the context of the \textit{Propagation}
function.  For a concise description of the parallel
\textit{Interpolation} function in what follows we denote by
\begin{equation} \label{eq:definitionneighbourboxesofpoint}
  \mathcal{M}^d(x) := \left\{ B^d_\mathbf{k} \in \mathcal{R}_B : x \in
    \mathcal{V} B^d_\mathbf{k} \right\},
\end{equation}
the set of \textit{cousin boxes of a surface discretization point
  $x \in \Gamma_N$ on level $d$}, $3 \leq d \leq D$, which extends the
concept of cousin boxes of a box, introduced in \eqref{eq:defneighboursandcousins} in
Section~\ref{sec:ifgfmethod}.
Using the definition \eqref{eq:definitionneighbourboxesofpoint}, the \textit{Parallel Interpolation} function is stated in Algorithm \ref{alg:ifgf_interpolation_parallel}. 
\begin{algorithm}[tb]
\small
\begin{algorithmic}[1]
  \ParFor{$x \in \Gamma_N$} \label{alg:parinterp_outerloop}
  \For{$B_\mathbf{k}^d \in \mathcal{M}^d(x)$}
  \State Determine $C_{\mathbf{k};\mathbf{\gamma}}^{d}$ s.t. $x \in C_{\mathbf{k};\mathbf{\gamma}}^{d}$ 
  \State Evaluate $I_P C_{\mathbf{k};\mathbf{\gamma}}^{d}(x) \times G(x, x_\mathbf{k}^d)$
  \EndFor 
  \EndParFor 
  \caption{Parallel Interpolation$(d)$}
\label{alg:ifgf_interpolation_parallel}
\end{algorithmic}
\end{algorithm}

In summary, the OpenMP parallelization strategies proposed above for
the functions \textit{Parallel LevelDEvaluations}, \textit{Parallel
  Propagation} and \textit{Parallel Interpolation} are thread-safe by
design, and they provide effective work distribution by relying on
iteration over items (relevant cone segments or surface discretization
points) that exist in a sufficiently large (and essentially constant)
quantities for all levels $d$, $3\leq d \leq D$, in the box-octree
structure. As a result, the proposed approach effectively eliminates
the hard limitation present in the straightforward OpenMP
parallelization scheme mentioned at the beginning of this
section. Note that the proposed IFGF box-cone parallelization strategy
is in general not applicable to other hierarchical acceleration
methods, such as e.g. FMM-type algorithms. Indeed, in contrast to the
incremental propagation and surface evaluation approach inherent in
the IFGF method, previous acceleration methods rely on the FFT
algorithm---which, as discussed in Section~\ref{sec:introduction},
leads to inefficiencies in the upper portions of the corresponding
octree structures~\cite{2007DirectionalFMMLexing,
  2014ParallelDirectionalFMMLexing}.

\subsection{IFGF MPI parallelization} \label{sec:mpiparallelization}
The proposed MPI parallel IFGF algorithm, which enables both data
distribution onto the MPI ranks and efficient communication of data
between MPI ranks, is described in detail in
Sections~\ref{sec:distribution} through~\ref{sec:communication}. The
approach mirrors the one proposed in
Section~\ref{sec:ompparallelization} for the corresponding OpenMP
interface. In fact, the MPI parallel scheme is based on slight
modifications of the OpenMP parallel
Algorithms~\ref{alg:ifgf_parallellevelD},
\ref{alg:ifgf_propagation_parallel}, and
\ref{alg:ifgf_interpolation_parallel}. As indicated by the theoretical
discussion in Section~\ref{sec:complexity}, the communication overhead
is such that the intrinsic IFGF linearithmic complexity previously
demonstrated in~\cite{Bauinger2021} for a single core implementation
is preserved on any fixed number $N_c$ of cores; an illustration of
this theoretical result on $N_c = 1,680$ cores is presented in the
Supplementary Materials Table~\ref{table:linearithmic}. Most
importantly, as in the OpenMP case (cf. the last paragraph of
Section~\ref{sec:ompparallelization}), for arbitrarily large numbers
$D$ of levels, the MPI IFGF algorithm iterates over items (relevant
cone segments or surface discretization points) that exist in a
sufficiently large (and essentially constant) quantities for all
levels $d$, $3\leq d \leq D$, in the box-octree structure. As a
result, the strategy results in an overall MPI-OpenMP IFGF parallel
scheme without hard limitations on the achievable parallelism as the
number of cores grows.


\subsubsection{Problem decomposition and data
  distribution} \label{sec:distribution}

The distribution of the data required by the IFGF algorithm to the MPI
ranks can be summarized as the independent distribution of the set of
surface discretization points $\Gamma_N$, which are organized on the
basis of boxes induced by an octree structure, and the distribution of
the set of relevant cone segments on each level
$\mathcal{R}_C^d$. Clearly, for an efficient parallel implementation,
the distribution used should balance the amount of work performed by
each rank while maintaining a minimal memory footprint per rank and
while also minimizing the communication between ranks. A concise
description of the method used for data distribution to the MPI ranks
is presented in what follows, where we let $N_r \in \mathbb{N}$ and
$\rho \in \mathbb{N}$ ($1 \leq \rho \leq N_r$) denote the number of
MPI ranks and the index of a specific MPI rank, respectively.

\begin{figure}
    \centering
    \includegraphics[width=.7\linewidth]{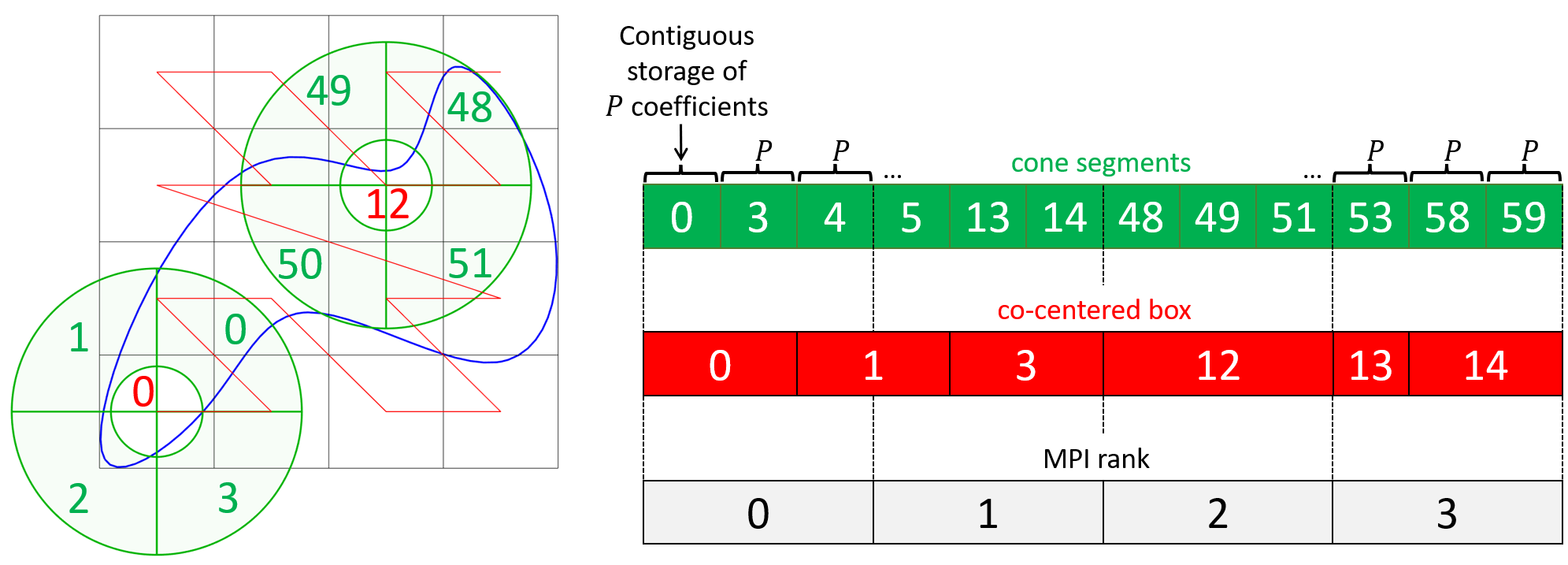}
    \caption{Left panel: Two-dimensional example of an ordering of the
      cone segments based on the Morton order of the boxes on level
      $d=3$ with four cone segments per box. The red line indicates
      the Morton order of the boxes where the red numbers denote the
      actual Morton code of the containing box. The green numbers
      denote the ordering of the cone segments in the proposed
      Morton-based cone-segment ordering. The blue curve denotes a
      sketch of a scatterer. Right panel: Sketch of a possible
      cone-segment memory layout, demonstrating the equi-distribution
      of cone segments among ranks, and emphasizing a central element
      of the proposed parallelization strategy, namely, that
      co-centered cone segments may be assigned to different MPI
      ranks.  Note that only relevant boxes and cone segments are
      stored in memory resulting in some numbers in the ordering being
      skipped.}
    \label{fig:sketchifgfordering}
\end{figure}

The distribution of the surface discretization points is orchestrated
on the basis of an ordering of the set of relevant boxes
$\mathcal{R}^d_B$ on each level $d$, which, in the proposed algorithm,
is obtained from a depth-first traversal of the octree structure. This
ordering is equivalent to a Morton order of the boxes (as described
e.g. in~\cite{2016Biros, 2013Warren, 2019Keyes,
  1993HashedOcttreeWarren} and depicted by the red
\reflectbox{Z}-looking curve in the left panel of
Figure~\ref{fig:sketchifgfordering}) which, as indicated in
\cite{1993HashedOcttreeWarren}, can be generated quickly from the
positions $\mathbf{k} \in I_B^D$ of the level-$D$ boxes
$B_\mathbf{k}^D$ through a bit-interleaving procedure. Ordering the
surface discretization points according to the Morton order of the
containing level-$D$ boxes also guarantees a Morton order on every
other level $d$, $1 \leq d \leq D-1$. More precisely, at every level
$d$ the Morton order introduces a total order $\prec$ on the set of
boxes. The ordering of the surface discretization points $\Gamma_N$ is
facilitated by assigning each point $x \in \Gamma_N$ the Morton order
of the containing level-$D$ box, which can be computed through a
division operation on the coordinates of the point $x$ to get the
index $\mathbf{k} \in I_B^D$ of the containing box with a subsequent
bit-interleaving procedure, followed by a simple sorting of the points
according to their assigned Morton order. Noting that the map which
assigns to each point on $\Gamma_N$ the Morton order of the containing
level-$D$ box is not injective, in order to obtain a total order on
all of $\Gamma_N$ we additionally order in an arbitrary manner subset
of points $x \in \Gamma_N$ with the same assigned Morton order.  The
resulting overall order has the desirable properties that, on every
level $d$, surface discretization points within any given box are
contiguous in memory, and that boxes close in real space are also
close in memory.

The sorted surface discretization points are distributed to the MPI ranks based on their containing level-$D$ boxes, in such a way that the boxes processed by each each rank are an ``interval'' set of the form 
$\{ B \in \mathcal{R}_B^D \ : \ B_\mathbf{k_1}^D \prec B \prec
B_\mathbf{k_2}^D \}$, for suitable choices of
$\mathbf{k_1}, \mathbf{k_2} \in I_B^D$ designed to guarantee that all the boxes on a given
rank contain a number of surface discretization points as close as
possible to the average value $N/N_r$. Hence, the smallest boxes in
the octree structure represent the smallest ``unit'' for the
distribution of the surface discretization points. The maximum
possible deviation in the number of discretization points assigned to
a certain MPI rank from the average is therefore given by the maximum
number of surface discretization points contained within one level-$D$
box in the octree structure. For reasonable distributions of the
discretization points $\Gamma_N$ on the surface $\Gamma$, and for a
suitable choice of the number of levels $D$, this deviation between
MPI ranks is typically less than 100 surface discretization points.

The set of \textit{surface discretization points stored in the
  $\rho$-th MPI rank}, $1 \leq \rho \leq N_r$, is denoted by
$\Gamma_{N,\rho}$. By definition, the subsets $\Gamma_{N,\rho}$ of
$\Gamma_N$ are pairwise disjoint and their union over all MPI ranks
$\rho = 1,\ldots,N_r$ equals $\Gamma_N$.  The distribution of the
surface discretization points is used to evenly divide between all MPI
ranks the work performed in the \textit{Interpolation} function
(OpenMP Algorithm~\ref{alg:ifgf_interpolation_parallel}).  The
underlying level-$D$ based distribution of $\Gamma_N$ is utilized
throughout all levels $d = D, \ldots, 3$. Thus, the MPI parallel
\textit{Interpolation} function results from the straightforward and
level-independent modification of Line \ref{alg:parinterp_outerloop}
in Algorithm \ref{alg:ifgf_interpolation_parallel}, to read
$x \in \Gamma_{N,\rho}$ instead of $x \in \Gamma_N$---as shown in
Algorithm \ref{alg:ifgf_interpolation_mpiparallel}. Naturally, the
values of the discrete operator $I(x_\ell)$ in~\eqref{eq:field1}
computed by the $\rho$-th MPI rank correspond to points
$x_\ell\in\Gamma_{N,\rho}$, and they are therefore also stored in the
$\rho$-th MPI rank. In other words, the set of resulting field values
$I(x_\ell)$ is partitioned and stored in the MPI ranks according to
the partition utilized for the surface discretization points
$\Gamma_N$.

\begin{algorithm}[tb]
\small
\begin{algorithmic}[1]
  \ParFor{$C_{\mathbf{k}; \mathbf{\gamma}}^D \in \mathcal{R}_{C, \rho}^D$}
  \State Evaluate and store $F_\mathbf{k}^D(\mathcal{X} C_{\mathbf{k};\mathbf{\gamma}}^D)$
  \State Generate interpolant $I_P C_{\mathbf{k}; \mathbf{\gamma}}^D$
  \EndParFor
  \caption{MPI Parallel LevelDEvaluations}
\label{alg:ifgf_mpiparallellevelD}
\end{algorithmic}
\end{algorithm}

\begin{algorithm}[tb]
\small
\begin{algorithmic}[1]
  \ParFor{$x \in \Gamma_{N,\rho}$}
  \For{$B_\mathbf{k}^d \in \mathcal{M}^d(x)$}
  \State Determine $C_{\mathbf{k};\mathbf{\gamma}}^{d}$ s.t. $x \in C_{\mathbf{k};\mathbf{\gamma}}^{d}$ 
  \State Evaluate $I_P C_{\mathbf{k};\mathbf{\gamma}}^{d}(x) \times G(x, x_\mathbf{k}^d)$
  \EndFor 
  \EndParFor 
  \caption{MPI Parallel Interpolation$(d)$}
\label{alg:ifgf_interpolation_mpiparallel}
\end{algorithmic}
\end{algorithm}

\begin{algorithm}[tb]
\small
\begin{algorithmic}[1]
  \ParFor{$C_{\mathbf{j}; \mathbf{\gamma}}^{d-1} \in \mathcal{R}_{C,\rho}^{d-1}$}
  \For{$B_\mathbf{k}^d \in \mathcal{C} (\mathcal{R}_B C_{\mathbf{j}; \mathbf{\gamma}}^{d-1})$}
  \For{$x \in \mathcal{X} C_{\mathbf{j};\mathbf{\gamma}}^{d-1}$}
  \State Determine $C_{\mathbf{k};\mathbf{\alpha}}^{d}$ s.t. $x \in C_{\mathbf{k};\mathbf{\alpha}}^{d}$ 
  \State Evaluate and add $I_P C_{\mathbf{k};\mathbf{\alpha}}^{d}(x) \times G(x, x_\mathbf{k}^d)/ G(x, x_\mathbf{j}^{d-1})$
  \EndFor 
  \EndFor 
  \State Generate interpolant $I_P C_{\mathbf{j}; \mathbf{\gamma}}^{d-1}$
  \EndParFor
  \caption{MPI Parallel Propagation$(d)$}
\label{alg:ifgf_propagation_mpiparallel}
\end{algorithmic}
\end{algorithm}

The data associated with the level-$d$ relevant cone segments is also
distributed to MPI ranks on the basis of a total order---in this case,
a total order on the set of level-$d$ cone segments that is based on
the Morton order imposed on the level-$d$ boxes, in such a way that
co-centered cone segments are close in memory. It should be noted
that, for every relevant cone segment
$C_{\mathbf{k};\mathbf{\gamma}}^d \in R_C^d$, $3\leq d \leq D$
(see~\eqref{eq:relconesegmentslevelD}), the set of $P$ coefficients
that characterize the polynomial interpolants
$I_P C_{\mathbf{k}; \bf{\gamma}}^d$ (Section~\ref{sec:ifgfmethod}),
which approximate the field $F_\mathbf{k}^d$ in~\eqref{eq:fieldperbox}
within the cone segment $C_{\mathbf{k}; \bf{\gamma}}^d$, need to be
stored, in appropriately distributed manner, for two consecutive
levels. Indeed, for each $d$, these level-$d$ coefficients are
utilized to enable two different interpolation procedures, namely
interpolation from level $d$ to interpolation points at the
parent-level $(d-1)$ in the \textit{Propagation} function (Line 4 in
Algorithm \ref{alg:ifgf_propagation_parallel}), as well as
interpolation to the level-$d$ cousin surface discretization points in
the \textit{Interpolation} function (Line 3 in Algorithm
\ref{alg:ifgf_interpolation_parallel}).


The set of level-$d$ relevant cone segments $\mathcal{R}^d_C$ is
sorted on the basis of the Morton order induced by the co-centered
level-$d$ boxes followed by a suitable sorting of cone segments in
each spherical coordinate system---resulting in a total order
$\sqsubset$ in the set of all level-$d$ relevant cone segments, as
depicted in the left panel of
Figure~\ref{fig:sketchifgfordering}. (Each set of co-centered cone
segments is ordered using the radial direction first, then elevation
and finally azimuth, although any other ordering could be used.)
Finally, at each level $d$ ($d = D, \ldots, 3$), approximately
equi-sized and pair-wise disjoint intervals of relevant cone segments
$C$ of the form
$\left\{C \in \mathcal{R}_C^d \ : \ C_{\mathbf{k_1};
    \bf{\gamma_1}}^d\sqsubset C \sqsubset C_{\mathbf{k_2};
    \bf{\gamma_2}}^d\right\}$, for some
$\mathbf{k_1}, \mathbf{k_2} \in I_B^d$ and
$\gamma_1, \gamma_2 \in I_C^d$ (i.e., disjoint intervals of {\em
  contiguous} cone segments), are distributed to the MPI ranks, as
illustrated in the right panel of
Figure~\ref{fig:sketchifgfordering}. Note that the specific assignment
of cone segments to MPI ranks is solely determined by the order
$\sqsubset$ and the number of MPI ranks and cone segments, and it does
not otherwise relate to the underlying box tree. In particular, as
suggested in the right panel of Figure~\ref{fig:sketchifgfordering},
co-centered cone segments may be assigned to different MPI
ranks---which induces a flexibility that leads to excellent
load-balancing and, therefore, high parallelization efficiency. As is
the case for the relevant boxes, the proposed ordering of the relevant
cone segments implies that cone segments which are close in real space
(i.e. co-centered with the same box and pointing in the same direction
or co-centered with boxes which are close in real space) are also
close in memory, and, in particular, are likely to be stored within
the same MPI rank. Analogously to the notation introduced above for
the distributed surface discretization points, the relevant level-$d$
cone segments assigned to a MPI rank with index $\rho$,
$1 \leq \rho \leq N_r$, are denoted by $\mathcal{R}_{C, \rho}^d$. The
MPI-capable algorithm is thus obtained by adjusting the loops in the
first lines in Algorithms \ref{alg:ifgf_parallellevelD} and
\ref{alg:ifgf_propagation_parallel} to only iterate over the level-$d$
relevant cone segments $\mathcal{R}_{C, \rho}^d$ stored in the current
rank $\rho$, as shown in the MPI parallel
Algorithms~\ref{alg:ifgf_mpiparallellevelD}
and~\ref{alg:ifgf_propagation_mpiparallel}, instead of iterating over
all relevant cone segments on level $d$.

\subsubsection{Practical implementation of the box-cone data
  structures\label{sec:datastructures}}
  
A C++ implementation of the parallel IFGF box-cone data structures described above, which enables a linearithmic memory and time complexity, is described in detail in what follows.\looseness = -1

In the proposed implementation the geometry $\Gamma_N$ is stored in
three separate arrays $X_1$, $X_2$ and $X_3$ (either C style arrays or
std::vector) of size $N$ for the $x_1$, $x_2$, and $x_3$ components of
the surface discretization points $(x_1, x_2, x_3) = x \in \Gamma_N$,
resulting in a \textit{structure of arrays} (SoA) memory layout
\cite{IntelSoA1}, which is beneficial as it leads to increased
floating-point performance under automatic vectorization on the basis
of \textit{single instruction, multiple data} (SIMD) hardware
\cite[Section 2.7]{sterling2017high} generally available in modern
processors. As mentioned above in Section~\ref{sec:distribution}, each one of the three arrays is sorted according to the
Morton order of the boxes. Similarly, the real and imaginary parts of the resulting field values $I(x_\ell)$, $1 \leq \ell \leq N$, are stored as two independent arrays, $I_\Re$ and $I_\Im$, of size $N$. The order of these field values coincides with the order imposed on the surface discretization points such that $I(x_\ell) = I_\Re[k] + \iota  I_\Im[k]$ at a given point $x_\ell$ is stored at the same position $k$ in the arrays $I_\Re$ and $I_\Im$ as the corresponding surface discretization point $x_\ell = (X_1[k], X_2[k], X_3[k])$ in the arrays $X_1$, $X_2$ and $X_3$.

The algorithm enacts the box-octree inherent in the IFGF solver in the form of a {\em
  linear octree structure}
(cf. \cite{cormen2009introduction,1990SpatialDataStructuresSamet}). In
particular, the proposed linear octree only includes
data associated with relevant boxes, and it does not store any
information about non-relevant boxes, to avoid $\mathcal{O}(N^{3/2})$ memory requirements, as described in detail in what
follows. Relevant boxes in the linear octree are represented,
on each level $d = 3, \ldots, D$, by the box index
$\mathbf{k} \in I_B^d$ (as described above in Section
\ref{sec:ifgfmethod}) and the equivalent Morton order. Each box stores the position in the arrays $X_1$, $X_2$, and $X_3$ of the first surface discretization point contained in the box in addition to the number of discretization points in the box in a hash map (cf. \cite[Section
11]{cormen2009introduction}) with average $\mathcal{O}(1)$ time and memory
complexity for read access (e.g. a std::unordered\_map), where the Morton order of the box is utilized as the key. Thus, given a Morton order of a box, the discretization points
contained within the box can be determined in
$\mathcal{O}(1)$ time and memory complexity. Conversely, given any
surface discretization point $x \in \Gamma_N$, the three-dimensional
index $\mathbf{k} \in I_B^d$ (for every level $d = 3, \ldots, D$) of
the box $B_\mathbf{k}^d$ containing the point $x$ and the associated Morton order can be determined
through simple division and bit-interleaving, as described in Section~\ref{sec:distribution}, in $\mathcal{O}(1)$ time and memory. 
Overall, this guarantees a true
$\mathcal{O}(N \log N)$ time and memory implementation by avoiding the
storage of any information regarding non-relevant boxes. Note that the
linear octree structure described above is essentially the same as the one
presented in \cite{1993HashedOcttreeWarren}.


Similarly, for each level $d$, the relevant cone segments
$C_{\mathbf{k}; \gamma}^d$, or, more precisely, the real and imaginary
parts of the $P$ coefficients representing the interpolants
$I_P C_{\mathbf{k}; \gamma}^d$ on the relevant cone segments
$C_{\mathbf{k}; \gamma}^d$, are stored separately in two
one-dimensional arrays per rank (following the above partition of the
cone segments to the ranks). To associate the three-dimensional cone
segment index $\gamma \in I_C^d$ with the actual coefficients, a hash
map for each relevant box $B_\mathbf{k}^d$ is used, where the
\textit{value} of the hash map is an index pointing to the first of
the coefficients $I_P C_{\mathbf{k}; \gamma}^d$ in the array of
coefficients mentioned above, and where $\gamma \in I_C^d$ is the
\textit{key} of the hash map. (Note that the three-dimensional cone
segment index $\gamma$, which runs over both relevant and non-relevant
cone segments, corresponds to the relative position of the cone
segment in the spherical coordinate system centered at the box
center.) The usage of a hash map circumvents the storage of any
non-relevant cone segment data while maintaining the association with
the three-dimensional index $\gamma$ that allows an easy
identification of the cone segment based on its relative position in
spherical coordinates. Thus, for a given Cartesian point
$x\in\mathbb{R}^3$, this data structure can be used to locate the
interpolant $I_P C_{\mathbf{k};\mathbf{\gamma}}^d$ for the relevant
cone segment $C_{\mathbf{k};\mathbf{\gamma}}^d \in R_C^d$ containing
the point $x$ through a transformation of $x$ to spherical coordinates
centered at the origin of the cone segment
$C_{\mathbf{k};\mathbf{\gamma}}^d$, a division to get the cone segment
index $\gamma$ and a look-up in the hash map to get the coefficients
of the interpolating polynomial. The association of any point with the
relevant cone segment containing it can therefore be achieved on
average in $\mathcal{O}(1)$ time and memory. This is required in the
\textit{Interpolation} and \textit{Propagation} function to facilitate
the interpolation to cousin surface discretization points and
parent-level interpolation points, respectively.

Note that, for increased performance, the hash maps described above and stored on a given rank $\rho$ are required to contain all associations between boxes, cone segments, discretization points and interpolation coefficients utilized by the current rank $\rho$ at any point in the algorithm. In particular, if certain surface discretization points or interpolant coefficients are stored on a different rank $\tilde \rho \neq \rho$, but are required in the current rank $\rho$, the above hash maps are utilized to find the data and, consequently, enable the communication of that data through MPI. While this produces some data duplication, analogously to ``halo regions'' \cite[Sec. 9.6]{sterling2017high} employed in grid-based methods, the memory duplicated in the parallel IFGF method is limited to surface discretization points and interpolant coefficients of neighbors and cousins.


\subsubsection{Data communication} \label{sec:communication}
Clearly, for an MPI rank to access data stored in
a different rank, explicit communication between the ranks must
take place. The proposed solution, which we favor due to the decreased
complexity of the implementation it entails, is based on
\textit{one-sided} or \textit{remote memory access} (RMA)
communication introduced in MPI-2 \cite[Section
5]{ParallelProgramming2013}, \cite[Section
8]{sterling2017high}---which utilizes a single MPI\_Get or MPI\_Put
call on the origin rank instead of a coupled MPI\_Recv-MPI\_Send call
(or similar functionalities) involving both the origin and the target
rank.

The data any MPI rank may require from other MPI ranks is limited to
certain interpolants $I_P C_{\mathbf{k};\mathbf{\gamma}}^d$. It is
therefore sufficient to store the corresponding coefficients
in so-called RMA \textit{windows} (in MPI given by
MPI\_Win and allocated with e.g. MPI\_Win\_allocate), which enable the
one-sided communication approach. For increased efficiency, the
computations and communications are organized among the ranks on the
basis of the following two considerations: 1)~For each $\rho$,
$1\leq \rho\leq N_r$, the $\rho$-th rank asynchronously collects from
other ranks all the data (i.e. the coefficients of the interpolants)
it requires to perform \textit{Interpolation} or \textit{Propagation}
computations assigned to it; and 2)~The communications necessary to
collect this data are interleaved with the computations in such a way
that while the computations by the \textit{Interpolation} function
take place, the communication for the next \textit{Propagation}
function is performed and vice versa.  This approach, which
effectively hides the communications behind computations (thus
increasing the performance and parallel efficiency), requires every
MPI rank to store all data it obtains from other ranks for one full
level-$d$ ($3\leq d \leq D$) \textit{Interpolation} or
\textit{Propagation} step while it continues to store the coefficients
it has itself generated---which effectively increases the peak memory
per rank requirements slightly (by e.g. 10$\%$ or less).

The level-$d$ dependent \textit{CommunicateInterpolationData} (resp.
\textit{CommunicatePropagationData}) programming function in
Algorithm~\ref{alg:ifgf_communicatepropagationdata}
(resp. Algorithm~\ref{alg:ifgf_communicateinterpolationdata})
encapsulates the communications performed by each rank to obtain, from
other ranks, the polynomial coefficients it needs to enact the
necessary level-$d$ interpolation computations (resp. interpolation
computations onto level-$(d-1)$ interpolation points) required by the
\textit{Interpolation} (resp. \textit{Propagation}) function. The
\textit{LevelDEvaluations} function does not need any communications
since the surface discretization points $x \in \Gamma_N$, which are required in the \textit{LevelDEvaluations} function but which are not stored as part of $\Gamma_N, \rho$ (see previous Section~\ref{sec:distribution}), are duplicated to the $\rho$-th MPI rank. The rank that stores a level-$D$ relevant cone
segment, as described at the beginning of Section
\ref{sec:mpiparallelization}, facilitates the evaluation of the field
at the interpolation points of that cone segment and the generation of
the interpolants independently from every other MPI rank.

Using the functions~\ref{alg:ifgf_mpiparallellevelD}
through~\ref{alg:ifgf_communicateinterpolationdata}, the pseudo-code
for the proposed overall MPI-OpenMP IFGF algorithm is given in
Algorithm~\ref{alg:ifgf_finalMPI}. Note that access to RMA windows is
usually asynchronous and requires some form of synchronization to
ensure the data transfer is finalized before the communicated data is
accessed. Moreover, the call to the
\textit{CommunicatePropagationData} in
Algorithm~\ref{alg:ifgf_finalMPI} requires for the
\textit{Propagation} function to have completed in all ranks targeted
by the communication function.

\begin{algorithm}[tb]
\small
\begin{algorithmic}[1]
  \ParFor{$C_{\mathbf{j}; \mathbf{\gamma}}^{d-1} \in \mathcal{R}_{C,\rho}^{d-1}$}
  \For{$B_\mathbf{k}^d \in \mathcal{C} (\mathcal{R}_B C_{\mathbf{j}; \mathbf{\gamma}}^{d-1})$}
  \For{$x \in \mathcal{X} C_{\mathbf{j};\mathbf{\gamma}}^{d-1}$}
  \State Find $\tilde \gamma$ such that $x \in C_{\mathbf{k};\mathbf{\tilde \gamma}}^{d} \in \mathcal{R}_C B_\mathbf{k}^d$
  \State Identify the MPI rank $\rho$ on which $I_P C_{\mathbf{k};\mathbf{\tilde \gamma}}^{d}$ is stored
  \State MPI\_Get $I_P C_{\mathbf{k};\mathbf{\tilde \gamma}}^{d}$ from rank $\rho$
  \EndFor 
  \EndFor 
  \EndParFor
  \caption{CommunicatePropagationData($d$)}
\label{alg:ifgf_communicatepropagationdata}
\end{algorithmic}
\end{algorithm}

\subsection{Parallel linearithmic complexity  analysis \label{sec:complexity}}
Reference~\cite[Sec. 3.3.3]{Bauinger2021} shows that the basic IFGF
algorithm runs on a linearithmic ($\mathcal{O}(N \log N)$) number of
arithmetic operations. The present section, in turn, shows that the
\textit{communication cost} additionally required by the proposed
MPI-OpenMP parallel IFGF algorithm also grows linearithmically---thus,
establishing that, on a fixed number of cores, the parallel algorithm
runs on an linearithmic overall computing time.

To do this, in view of the data distribution strategy described in
Section~\ref{sec:distribution}, it suffices to ensure that both the
\textit{Interpolation} and \textit{Propagation} functions require a
linearithmic communication cost. Inspection of the corresponding
Algorithms~\ref{alg:ifgf_interpolation_mpiparallel}
and~\ref{alg:ifgf_propagation_mpiparallel} (specifically, lines 4 and
5, respectively) shows that these functions, and, thus, the overall
parallel IFGF algorithm, only require communication of certain
polynomial coefficients---a task that is effected via the
communication Algorithms~\ref{alg:ifgf_communicateinterpolationdata}
and~\ref{alg:ifgf_communicatepropagationdata}, respectively. Thus, the
analysis of the communication cost amounts to counting the number of
coefficients that are communicated, including multiple counts for
coefficients that are communicated to multiple ranks, as a result of
the application of these two communication algorithms within the
overall IFGF algorithm.

In order to count the number of communications effected by each one of
these algorithms we proceed as follows.  Noting that, since, 1)~As
indicated in~\cite[Sec. 3.3.3]{Bauinger2021}, there are
$\mathcal{O}(N)$ relevant cone segments per level, each one of which
contains $\mathcal{O}(1)$ data (namely, the $P$ coefficients of a
single polynomial interpolant); 2)~Each cone-segment data is stored in
exactly one MPI rank (Section~\ref{sec:distribution}); and, as
discussed below for both communication algorithms, 3)~Each relevant
cone segment is communicated to a uniformly bounded number of MPI
ranks at each level $d = 3, \ldots, D$; it follows that for each level
$d$ ($3\leq d\leq D$) a total of $\mathcal{O}(N)$ coefficients are
communicated by each of the communication
algorithms~\ref{alg:ifgf_communicatepropagationdata}
and~\ref{alg:ifgf_communicateinterpolationdata} for each one of the
$D=\log{N}$ levels, at a total communication cost of
$\mathcal{O}(N\log N)$ coefficients by this algorithm, as desired.

\begin{algorithm}[tb]
\small
\begin{algorithmic}[1]
  \ParFor{$x \in \Gamma_{N, \rho}$}
  \For{$B_\mathbf{k}^d \in \mathcal{M}^d(x)$}
  \State Find $\tilde \gamma$ such that $x \in C_{\mathbf{k};\mathbf{\tilde \gamma}}^{d} \in \mathcal{R}_C B_\mathbf{k}^d$
  \State Identify the MPI rank $\rho$ on which $I_P C_{\mathbf{k};\mathbf{\tilde \gamma}}^{d}$ is stored
  \State MPI\_Get $I_P C_{\mathbf{k};\mathbf{\tilde \gamma}}^{d}$ from rank $\rho$
  \EndFor 
  \EndParFor 
  \caption{CommunicateInterpolationData($d$)}
\label{alg:ifgf_communicateinterpolationdata}
\end{algorithmic}
\end{algorithm}

It remains for us to show that point~3) above holds for both
communication algorithms. In the case of the propagation communication
we note that each relevant cone segment on any level
$d = D, \ldots, 4$ is split into eight smaller cone segments on the
parent level $(d-1)$. Thus, for each level-$d$ relevant cone segment,
this results in at most $K$ parent-level cone segments (usually $K=8$,
or possibly a slightly higher number owing to the re-centering
procedure associated with the \textit{Propagation} function, but most
often $K=1$) that could be targets for the interpolation procedure in
the \textit{Propagation} function. In view of point~2) above, each
level-$d$ relevant cone segment must thus communicate coefficients to
no more than $\mathcal{O}(1)$ ranks, and point~3) follows in this
case.

In the case of the interpolation communication, finally, relevant
cone-segment coefficients need to be communicated to ranks that store
surface discretization points included in boxes that are cousins of
the box co-centered with the relevant cone segment. First, on the
lowest level $D$, each relevant box has at most $K = 189$ cousin boxes
and since, by design, the surface discretization points contained
within each one of the smallest boxes are stored in a single MPI rank
(Section~\ref{sec:distribution}), it follows that $\mathcal{O}(1)$ (at
most $189$) different MPI ranks require coefficients contained in each
relevant cone segment. Further, since each cone segment is partitioned
into eight in the transition from a given level $d$ to a subsequent
level $(d-1)$ (so that the number of relevant level-$D$ boxes contained
within a level-$(d-1)$ cone equals approximately one-fourth of the
corresponding number for level-$d$ cone segments, since $\Gamma_N$ is
a discretization of a 2D surface), and since, conversely, the number
of MPI ranks storing surface discretization points within a cousin box
increases by approximately a factor four in the same $d$-to-$(d-1)$
transition, the number of communications per relevant cone segment
remains essentially constant as a result of the $d$-to-$(d-1)$ level
transition. It follows that each relevant cone segment is communicated
to a $\mathcal{O}(1)$ number of MPI ranks for all levels $d$, thus
establishing the validity of point~3) for the interpolation
communication function, and completing the proof of linearithmic
complexity of the proposed parallel IFGF algorithm.

\begin{algorithm}[tb]
\small
\begin{algorithmic}[1]
  \State LevelDEvaluations()
  \State CommunicatePropagationData($D$) 
  \State
  \For{$d = D, \ldots, 3$}
  \State CommunicateInterpolationData($d$)
  \If {$d > 3$}
  \State Propagation($d$)
  \If {$d > 4$}
  \State CommunicatePropagationData($d-1$) 
  \EndIf
  \EndIf 
  \State Interpolation($d$) 
  \EndFor
  \caption{IFGF Method}
\label{alg:ifgf_finalMPI}
\end{algorithmic}
\end{algorithm}

\section{Numerical Results\label{sec:results}}
\begin{figure}
  \centering
  \includegraphics[width=.43\linewidth]{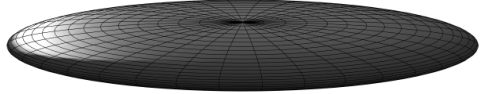}\hspace{1cm}\includegraphics[width=0.35\textwidth]{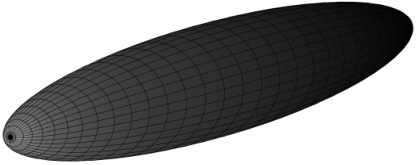}
  \caption{Test geometries. Left: Oblate spheroid
    $x^2 + y^2 + (z/0.1)^2 = a^2$. Right: Prolate spheroid
    $x^2 + y^2 + (z/10)^2 = a^2$.\label{fig:geometries}}
\end{figure}

Our numerical examples focus on three simple geometries which coincide
with the test cases presented in \cite{Bauinger2021}: a sphere of
radius $a$, the oblate spheroid $x^2 + y^2 + (z/0.1)^2 = a^2$ and the
prolate spheroid $x^2 + y^2 + (z/10)^2 = a^2$. The latter two
geometries are depicted in Figure~\ref{fig:geometries}. In what
follows the diameter (also referred to as the ``size'') of a geometry
$\Gamma$ is denoted by
\begin{equation} \label{eq:largestdiameter}
    d := d(\Gamma) := \max \limits_{x, y \in \Gamma} |x - y|,
\end{equation}
(not to be confused with the level index $d$ introduced in
Section~\ref{sec:ifgfmethod}); clearly we have $d = 2a$ in the case of
the sphere and the oblate spheroid geometries and $d = 20a$ for the
prolate spheroid geometry.  These relatively simple geometries present
the same kinds of challenges, in the context of the IFGF method, that
arise in a wide range of real-world problems, including aircraft,
lenses and meta-materials (with a point distribution somewhat similar
to that in an oblate spheroid), submarines (prolate spheroid), etc. For
example, even though the problem of finding a scattering solution for
a submarine is much more challenging than the corresponding problem for a
spheroid of the same size, in view of the need for accurate
integration of singular kernels and adequate representation of the
surface Jacobians, the performance of the IFGF method for the
evaluation of the discrete operator~\eqref{eq:field1} for a submarine
should not differ significantly from the corresponding performance on
a prolate spheroid of a comparable discretization, point distribution
and electromagnetic size.

In what follows we present IFGF performance data based on various runs
for these geometries. For our examples we utilize discretizations
$\Gamma_N$ obtained from use of parametrized surface patches covering
$\Gamma$ and equispaced partitioning of the corresponding parameter
spaces. The computer used is described briefly prior to the beginning
of Section~\ref{sec:ompparallelization} and, in more detail, in the
Supplementary Materials Section~\ref{sec:hardware}. The strong and
weak efficiency and speedup scalability concepts are detailed in
Section~\ref{sec:defeff}; briefly, relative to a base core-number
$N_c^0$, the $N_c$-core run speedup $S_{N_c^0, N_c}$ and the weak and
strong efficiencies $E^w_{N_c^0, N_c}$ and $E^s_{N_c^0, N_c}$ are used
to characterize the effectiveness of the proposed parallelization
schemes by relating computing times and core numbers under
weak-scaling tests (in which $N_c$ is increased proportionally to the
size $N$ of the discretization $\Gamma_N$) and strong-scaling tests
(wherein $N_c$ is increased as $N$ is held fixed).

\begin{table} [h]
  \centering
  \small
\begin{tabular}{|| c | r | c | c | c | l | l | r | r ||} 
 \hline
\multicolumn{1}{|| c|}{$\Gamma$} & \multicolumn{1}{| c |}{\bf $N$} & \multicolumn{1}{|c|}{$d$} & \multicolumn{1}{|c|}{ Nodes} & \multicolumn{1}{|c|}{\bf $N_c$} & \multicolumn{1}{|c|}{$\mathbf \varepsilon$} & \multicolumn{1}{|c|}{ $T$ (s)} & \multicolumn{1}{|c|}{ $E_{56, N_c}^w$ (\%)} & \multicolumn{1}{|c|}{\bf $E_{\frac{N_c}{4}, N_c}^w$ } \Tstrut \Bstrut \\ \hline \hline
 \multirow{3}{*}{Sphere} & $1,572,864$ &$128 \lambda$ & $1$ & $56$ &  $2 \times 10^{-3}$ & $7.77\times 10^{1}$ & $100$ & - \Tstrut \\
  &$6,291,456$ & $256 \lambda$ & $4$ & $224$ & $2\times 10^{-3}$ &  $9.78\times 10^{1}$ & $87$ & $87$ \Strut \\
  &$25,165,824$ & $512 \lambda$ & $16$ & $896$ & $2\times 10^{-3}$ &  $1.34\times 10^{2}$ & $69$ & $79$ \Bstrut \\
 \hline
 \multirow{3}{*}{\shortstack{Oblate\\Spheroid}} & $1,572,864$ &$128 \lambda$ & $1$ & $56$ &  $7\times 10^{-4}$ & $2.99\times 10^{1}$ & $100$ & - \Tstrut \\
  &$6,291,456$ & $256 \lambda$ & $4$ & $224$ & $6\times 10^{-4}$ &  $4.17\times 10^{1}$ & $79$ & $79$ \Strut \\
  &$25,165,824$ & $512 \lambda$ & $16$ & $896$ & $8\times 10^{-4}$ &  $5.74\times 10^{1}$ & $62$ & $79$ \Bstrut \\
 \hline
 \multirow{3}{*}{\shortstack{Prolate\\Spheroid}} & $6,291,456$ &$256 \lambda$ & $1$ & $56$ &  $5\times 10^{-4}$ & $4.97\times 10^{1}$ & $100$ & - \Tstrut \\
  &$25,165,824$ & $512 \lambda$ & $4$ & $224$ & $6\times 10^{-4}$ &  $6.83\times 10^{1}$ & $79$ & $79$ \Strut \\
  &$100,663,296$ & $1,024 \lambda$ & $16$ & $896$ & $7\times 10^{-4}$ &  $9.29\times 10^{1}$ & $63$ & $79$ \Bstrut \\
 \hline
\end{tabular}
\caption{Weak scaling test transitioning from 1 to 4 nodes, and then
  from 4 to 16 nodes, for three different geometries. The number of
  nodes, each one contaning $N_c=56 $ cores, is kept proportional to
  the number of surface discretization points, as required by the
  weak-scaling paradigm. }
\label{table:weakmpi}
\end{table}
\normalsize

Table~\ref{table:weakmpi} demonstrates the weak IFGF parallel
efficiency, for all three geometries considered, from a single compute
node ($N_c^0=56$) to 4 and 16 compute nodes ($N_c = 224$ and $896$,
respectively). We find that the efficiency relative to the base
$N_c^0=56$ case steadily decreases, but, importantly, the weak
relative efficiency $E_{\frac{N_c}{4}, N_c}^w$ remains essentially
constant as $N_c$ increases. Thus, under the assumption that this
trend is maintained for arbitrarily large numbers of nodes (as is
expected in view of the discussion in the first paragraph of
Section~\ref{sec:mpiparallelization} concerning absence of hard
limitations on achievable parallelism), the parallel IFGF method is
applicable to arbitrarily large problems---provided correspondingly
large hardware is used---with a constant $\approx 80\%$ efficiency
factor (cf. Table~\ref{table:weakmpi}) as the problem and hardware
sizes are both quadrupled from a given point of
reference. Section~\ref{sec:strong} demonstrates a similar quality of
the proposed algorithm under strong-scaling tests.

\begin{figure}
    \centering
    \includegraphics[width=.6\textwidth]{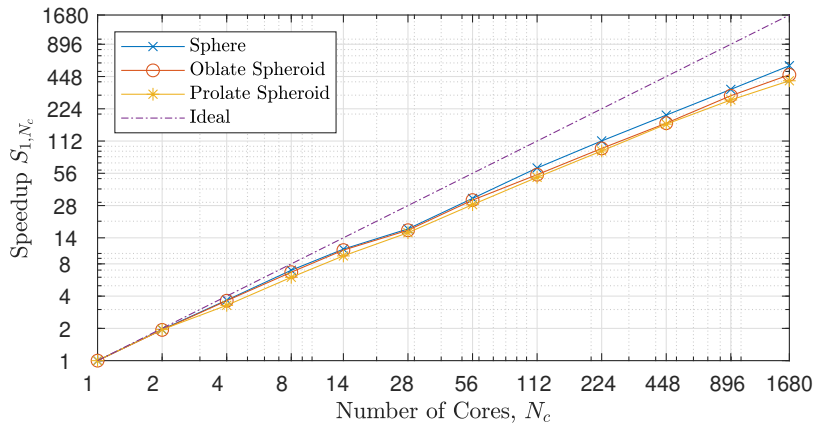}
    \caption{Measured speedup $S_{1, N_c}$ (vertical axis) versus
      number of cores $N_c$ (horizontal axis) in a strong scaling test
      transitioning from 1 core to 1,680 cores (= 30 compute nodes)
      for three geometries: a sphere of size 128 wavelengths (blue),
      an oblate spheroid of size 128 wavelengths (red), and prolate
      spheroid of size 256 wavelengths (yellow). The dash-dotted
      purple line indicates the theoretical perfect speedup. }
    \label{fig:strongscalingspeedupall}
  \end{figure}
The observed speedups under strong scaling tests, in turn, are
displayed in Figure~\ref{fig:strongscalingspeedupall}; additional
details concerning these results are provided in
Section~\ref{sec:strong}.  This figure presents speedup tests for
three test cases: a sphere of diameter $d = 128 \lambda$ (where
$\lambda =\frac{2\pi}{\kappa}$ denotes the wavelength and $d$ is given
in \eqref{eq:largestdiameter}), and oblate and prolate spheroids
(Figure~\ref{fig:geometries}) of large diameters $d = 128\lambda$ and
$d = 256 \lambda$, respectively.  The curves in
Figure~\ref{fig:strongscalingspeedupall} display, in each case, the
observed speedup $S_{1, N_c}$ for $1\leq N_c \leq 1,680$.  In view of
the requirements of the strong-scaling setup, test problems were
selected that can be run in a reasonable time on a single core and
with the memory available in the corresponding compute node. Clearly,
such test problems tend to be too small to admit a perfect
distribution onto large numbers of cores.  As illustrated in
Figure~\ref{fig:strongscalingspeedupall}, however, in spite of this
constraint, excellent scaling is observed in the complete range going
from 1 core to $1,680$ cores (30 nodes). As in the weak-scaling tests,
further, there is no hard limitation on scaling, even for such small
problems, (once again, in line with the discussion presented in
Section~\ref{sec:parallelization}), and it is reasonable to expect
that, unlike other approaches (for which either hard limits
arise~\cite{2014ParallelDirectionalFMMLexing} as described in the
first paragraph of Section~\ref{sec:ompparallelization}, or which rely on
memory duplication~\cite{2019LargeMLFMA,2021LargeMLFMA}), the
observed speedup continues to scale with the number of cores, as
suggested by Figure~\ref{fig:strongscalingspeedupall}, up to very
large numbers of cores. The computing speedups achieved by the
proposed parallel strategy outperform those achieved by other
MPI-parallel implementations of FMM and other numerical
methods~\cite{2009LargeSingleLvlFMM, 2009Ergul, 2020LargeFDTD}, and can be best
appreciated by noting that, instead of the e.g. approximately 40
minutes ($2.54\times 10^3$ secs., see first line in
Table~\ref{table:strongomp128} in Section~\ref{sec:strong}) required
by a single-core IFGF run, a total of $4.5$ secs.
($4.5 = 2.54\times 10^3/S_{1,1680}$ secs., where, per
Figure~\ref{fig:strongscalingspeedupall}, $S_{1,1680} = 565$) suffices
for the corresponding 1,680-core IFGF run.  It is interesting to note
that an approximately $1.51$ second 1,680-core run would have resulted
under perfect scaling.

\begin{figure}
    \centering
    \includegraphics[width=.6\textwidth]{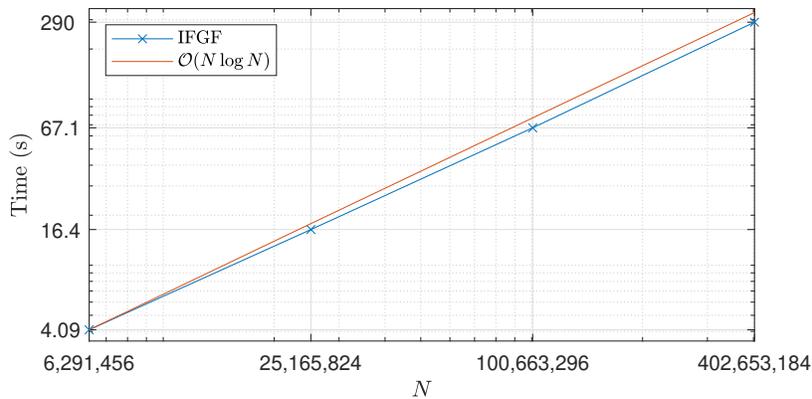}
    \caption{Illustration of the linearithmic complexity of the
      parallel IFGF method (which had previously been
      demonstrated~\cite{Bauinger2021} for the serial version of the
      algorithm), for the prolate spheroid geometry, on 30 compute
      nodes, with error $\varepsilon\approx 1.5\times 10^{-2}$. The
      acoustic diameter of the ellipsoid is kept proportional to
      $\sqrt{N}$, and it ranges from $512\lambda$ to
      $4,096\lambda$. Clearly, the parallel implementation preserves
      (and, in fact, slightly improves upon) the ideal linearithmic
      scaling. For this test one MPI rank per node and 56 OpenMP
      threads per MPI rank were used (resulting in 1680 cores). The
      peak IFGF memory used per MPI rank (excluding the memory
      required to store the initial geometry) as well as other
      additional data in tabular form are presented in
      Table~\ref{table:linearithmic}. }
\label{fig:linearithmic}
\end{figure}
Concluding this section, Figure~\ref{fig:linearithmic} presents
results of an investigation regarding the linearithmic scaling of the
parallel IFGF method for the prolate spheroid geometry on a fixed
number of nodes, namely, all 30 nodes available in the computer
cluster we use, and for $N$ ranging from 6,291,456 to 402,653,184, for
corresponding diameters ranging from $512\lambda$ to $4,096\lambda$.
The data in this figure, which is also presented in tabular form in
Table~\ref{table:linearithmic}, was generated by pinning a single MPI
rank to each compute node, each of which spawns 56 OpenMP threads,
with parameters resulting in an IFGF error
$\varepsilon\approx 1.5\times 10^{-2}$
(cf. equation~\eqref{eq:errorcomputation_2}).  The results show that
the linearithmic algorithmic complexity and memory requirements of the
basic IFGF algorithm are maintained in the parallel setting. Indeed,
the observed complexity even slightly outperforms the postulated
$\mathcal{O}(N \log N)$ within this range of values of $N$;
cf. Table~\ref{table:linearithmic} which suggests convergence to exact
linearithmic complexity as $N$ grows.

\section{Concluding Remarks} \label{sec:conclusio} This paper
presented a parallel version of the IFGF acceleration method
introduced in~\cite{Bauinger2021}, demonstrating in practice excellent
parallel scaling to large core numbers while simultaneously preserving
the linearithmic complexity of the sequential IFGF algorithm. The
proposed parallelization approach exploits the box-cone octree
structure inherent in the IFGF method, resulting in a strategy that,
per the theoretical discussion in Section~\ref{sec:ompparallelization}
and in the first paragraph of Section~\ref{sec:mpiparallelization}, is
applicable to arbitrarily large number of processing cores, and it
thereby does not suffer from bottlenecks or hard limits inherent in
approaches that orchestrate the parallelization on the basis of
octree-box partitioning only. A number of additional questions are
left for future work, as briefly mentioned in what follows. On one
hand, the feasibility of implementations on heterogeneous
architectures such as, e.g., computer systems that incorporate
general purpose graphical processing units (GPUs), is currently under
study. In particular, the use of GPUs to accelerate the interpolation
processes, which represent the most time consuming part of the IFGF
method, appears as highly promising avenue of inquiry. Additionally,
minor modifications to the data-decomposition strategy introduced in
Section~\ref{sec:distribution} could be introduced to not only
(approximately) equipartition the surface discretization points and cone segments among
MPI ranks, but also incorporate the number of actual computations
and the amount of data required from other MPI ranks in the partitioning scheme. Such an improved
data-decomposition design could indeed be obtained by relying on minor
adjustments to the cone and box intervals introduced in
Section~\ref{sec:distribution} leading to improved load-balancing,
and, thus, improved parallel efficiency.

\section*{Acknowledgments}
This work was supported by NSF, DARPA and AFOSR through contracts
DMS-2109831 and HR00111720035, FA9550-19-1-0173 and FA9550-21-1-0373,
and by the NSSEFF Vannevar Bush Fellowship under contract number
N00014-16-1-2808.

\bibliographystyle{unsrt}
\bibliography{IFGF}

\begin{thebibliography}{10}

\bibitem{Bauinger2021}
Christoph Bauinger and Oscar Bruno.
\newblock ``{I}nterpolated {F}actored {G}reen {F}unction'' method for
  accelerated solution of scattering problems.
\newblock {\em Journal of Computational Physics}, 430, 01 2021.

\bibitem{2007DirectionalFMMLexing}
Björn Engquist and Lexing Ying.
\newblock Fast directional multilevel algorithms for oscillatory kernels.
\newblock {\em SIAM Journal on Scientific Computing}, 29(4):1710--1737, 2007.

\bibitem{2012FMMChebyshevMessnerSchanz}
Matthias Messner, Martin Schanz, and Eric Darve.
\newblock Fast directional multilevel summation for oscillatory kernels based
  on chebyshev interpolation.
\newblock {\em Journal of Computational Physics}, 231:1175--1196, 2012.

\bibitem{Cands2009FastButterfly}
Emmanuel~J. Cand{\`e}s, Laurent Demanet, and Lexing Ying.
\newblock A fast butterfly algorithm for the computation of fourier integral
  operators.
\newblock {\em Multiscale Model. Simul.}, 7:1727--1750, 2009.

\bibitem{1996ButterflyMichielssen}
Eric Michielssen and Amir Boag.
\newblock A multilevel matrix decomposition algorithm for analyzing scattering
  from large structures.
\newblock {\em IEEE Transactions on Antennas and Propagation},
  44(8):1086--1093, 1996.

\bibitem{ROKHLIN1993}
V.~Rokhlin.
\newblock Diagonal forms of translation operators for the helmholtz equation in
  three dimensions.
\newblock {\em Applied and Computational Harmonic Analysis}, 1(1):82 -- 93,
  1993.

\bibitem{2006FMMRokhlin}
Hongwei Cheng, William~Y. Crutchfield, Zydrunas Gimbutas, Leslie~F. Greengard,
  J.~Frank Ethridge, Jingfang Huang, Vladimir Rokhlin, Norman Yarvin, and
  Junsheng Zhao.
\newblock A wideband fast multipole method for the helmholtz equation in three
  dimensions.
\newblock {\em Journal of Computational Physics}, 216:300--325, 2006.

\bibitem{2017Boerm}
Steffen B\"{o}rm and Jens Melenk.
\newblock Approximation of the high-frequency helmholtz kernel by nested
  directional interpolation.
\newblock {\em Numerische Mathematik}, 137(1):1--37, 10 2017.

\bibitem{2017BoermH2}
Steffen B\"{o}rm.
\newblock Directional h2‐matrix compression for high‐frequency problems.
\newblock {\em Numerical Linear Algebra with Applications}, 24, 07 2017.

\bibitem{2001FFTKunyansky}
Oscar~P. Bruno and Leonid~A. Kunyansky.
\newblock A fast, high-order algorithm for the solution of surface scattering
  problems: Basic implementation, tests, and applications.
\newblock {\em Journal of Computational Physics}, 169:80--110, 2001.

\bibitem{1996AIMFFTBelszynski}
E.~{Bleszynski}, M.~{Bleszynski}, and T.~{Jaroszewicz}.
\newblock Aim: Adaptive integral method for solving large-scale electromagnetic
  scattering and radiation problems.
\newblock {\em Radio Science}, 31(5):1225--1251, 1996.

\bibitem{1997FFTPhillips}
Joel~R. Phillips and Jacob~K. White.
\newblock A precorrected-fft method for electrostatic analysis of complicated
  3-d structures.
\newblock {\em IEEE Transactions on computer-aided design of integrated
  circuits and systems}, 16(10):1059--1072, 1997.

\bibitem{2003FMMLexingKernelIndependent}
Lexing Ying, George Biros, Denis Zorin, and M.~Harper Langston.
\newblock A new parallel kernel-independent fast multipole method.
\newblock In {\em A New Parallel Kernel-Independent Fast Multipole Method}, 11
  2003.

\bibitem{2003BebendorfRjasanow}
Mario Bebendorf and Sergej Rjasanow.
\newblock Adaptive low-rank approximation of collocation matrices.
\newblock {\em Computing}, 70:1--24, 02 2003.

\bibitem{2014ParallelDirectionalFMMLexing}
Austin~R. Benson, Jack Poulson, Kenneth Tran, Björn Engquist, and Lexing Ying.
\newblock A parallel directional fast multipole method.
\newblock {\em SIAM Journal on Scientific Computing}, 36(4):C335--C352, 2014.

\bibitem{ParallelFMMChandramowlishwaran2010}
A.~{Chandramowlishwaran}, S.~{Williams}, L.~{Oliker}, I.~{Lashuk}, G.~{Biros},
  and R.~{Vuduc}.
\newblock Optimizing and tuning the fast multipole method for state-of-the-art
  multicore architectures.
\newblock In {\em 2010 IEEE International Symposium on Parallel Distributed
  Processing (IPDPS)}, pages 1--12, 2010.

\bibitem{2021LargeMLFMA}
Rui-Qing Liu, Xiao-Wei Huang, Yu-Lin Du, Ming-Lin Yang, and Xin-Qing Sheng.
\newblock Massively parallel discontinuous galerkin surface integral equation
  method for solving large-scale electromagnetic scattering problems.
\newblock {\em IEEE Transactions on Antennas and Propagation},
  69(9):6122--6127, 2021.

\bibitem{2019LargeMLFMA}
Ming-Lin Yang, Yu-Lin Du, and Xin-Qing Sheng.
\newblock Solving electromagnetic scattering problems with over 10 billion
  unknowns with the parallel mlfma.
\newblock In {\em 2019 Photonics Electromagnetics Research Symposium - Fall
  (PIERS - Fall)}, pages 355--360, 2019.

\bibitem{2009Ergul}
{\"O}zg{\"u}r Erg{\"u}l and Levent Gurel.
\newblock A hierarchical partitioning strategy for an efficient parallelization
  of the multilevel fast multipole algorithm.
\newblock {\em IEEE Transactions on Antennas and Propagation},
  57(6):1740--1750, 2009.

\bibitem{2007Volakis}
Caleb Waltz, Kubilay Sertel, Michael~A. Carr, Brian~C. Usner, and John~L.
  Volakis.
\newblock Massively parallel fast multipole method solutions of large
  electromagnetic scattering problems.
\newblock {\em IEEE Transactions on Antennas and Propagation},
  55(6):1810--1816, 2007.

\bibitem{2009LargeSingleLvlFMM}
Luis Landesa, Jose Taboada, Fernando Obelleiro, Jose Rodriguez, José~Carlos
  Gallego, and Andres Gomez.
\newblock Solution of very large integral‐equation problems with
  single‐level fmm.
\newblock {\em Microwave and Optical Technology Letters}, 51:2451 -- 2453, 10
  2009.

\bibitem{2014Fanghzou1}
Fangzhou Wei and Ali~E. Yılmaz.
\newblock A more scalable and efficient parallelization of the adaptive
  integral method—part i: Algorithm.
\newblock {\em IEEE Transactions on Antennas and Propagation}, 62(2):714--726,
  2014.

\bibitem{2014Fangzhou2}
Fangzhou Wei and Ali~E. Yılmaz.
\newblock A more scalable and efficient parallelization of the adaptive
  integral method—part ii: Bioem application.
\newblock {\em IEEE Transactions on Antennas and Propagation}, 62(2):727--738,
  2014.

\bibitem{Nikolic2014}
Milo{\v{s}} Nikoli{\'{c}}, Aleksandar Jovi{\'{c}}, Josip Jaki{\'{c}}, Vladimir
  Slavni{\'{c}}, and Antun Bala{\v{z}}.
\newblock {\em An Analysis of FFTW and FFTE Performance}, pages 163--170.
\newblock Springer International Publishing, Cham, 2014.

\bibitem{Gumerov2004}
Nail~A. Gumerov and Ramani Duraiswami.
\newblock {\em Fast Multipole Methods for the Helmholtz Equation in Three
  Dimensions}.
\newblock Elsevier Science, 2004.

\bibitem{2019Keyes}
Mustafa Abduljabbar, Mohammed~Al Farhan, Noha Al-Harthi, Rui Chen, Rio Yokota,
  Hakan Bagci, and David Keyes.
\newblock Extreme scale fmm-accelerated boundary integral equation solver for
  wave scattering.
\newblock {\em SIAM Journal on Scientific Computing}, 41(3):C245--C268, 2019.

\bibitem{2014ButterflyLexing}
Jack Poulson, Laurent Demanet, Nicholas Maxwell, and Lexing Ying.
\newblock A parallel butterfly algorithm.
\newblock {\em SIAM Journal on Scientific Computing}, 36(1):C49--C65, 2014.

\bibitem{2012Winkel}
Mathias Winkel, Robert Speck, Helge H{\"u}bner, Lukas Arnold, Rolf Krause, and
  Paul Gibbon.
\newblock A massively parallel, multi-disciplinary barnes--hut tree code for
  extreme-scale n-body simulations.
\newblock {\em Computer physics communications}, 183(4):880--889, 2012.

\bibitem{2016Biros}
Dhairya Malhotra and George Biros.
\newblock Algorithm 967: A distributed-memory fast multipole method for volume
  potentials.
\newblock {\em ACM Transactions on Mathematical Software (TOMS)}, 43(2):1--27,
  2016.

\bibitem{2013Warren}
Michael~S. Warren.
\newblock 2hot: An improved parallel hashed oct-tree n-body algorithm for
  cosmological simulation.
\newblock In {\em SC '13: Proceedings of the International Conference on High
  Performance Computing, Networking, Storage and Analysis}, pages 1--12, 2013.

\bibitem{jimenez2021ifgf}
Edwin Jimenez, Christoph Bauinger, and Oscar~P Bruno.
\newblock Ifgf-accelerated integral equation solvers for acoustic scattering.
\newblock {\em arXiv preprint arXiv:2112.06316}, 2021.

\bibitem{Bull2010}
Mark Bull, James Enright, Xu~Guo, Chris Maynard, and Fiona Reid.
\newblock Performance evaluation of mixed-mode {OpenMP/MPI} implementations.
\newblock {\em International Journal of Parallel Programming}, 38:396--417, 10
  2010.

\bibitem{Drosinos2004}
N.~Drosinos and N.~Koziris.
\newblock Performance comparison of pure mpi vs hybrid mpi-openmp
  parallelization models on smp clusters.
\newblock In {\em 18th International Parallel and Distributed Processing
  Symposium, 2004. Proceedings.}, pages 15--, 2004.

\bibitem{Akhmetova2007}
Dana Akhmetova, Roman Iakymchuk, Orjan Ekeberg, and Erwin Laure.
\newblock Performance study of multithreaded mpi and openmp tasking in a large
  scientific code.
\newblock In {\em 2017 IEEE International Parallel and Distributed Processing
  Symposium Workshops (IPDPSW)}, pages 756--765, 2017.

\bibitem{dekking2005modern}
F.M. Dekking, C.~Kraaikamp, H.P. Lopuha{\"a}, and L.E. Meester.
\newblock {\em A Modern Introduction to Probability and Statistics:
  Understanding Why and How}.
\newblock Springer Texts in Statistics. Springer, 2005.

\bibitem{reinders2021data}
James Reinders, Ben Ashbaugh, James Brodman, Michael Kinsner, John Pennycook,
  and Xinmin Tian.
\newblock {\em Data Parallel C++: Mastering DPC++ for Programming of
  Heterogeneous Systems using C++ and SYCL}.
\newblock Springer Nature, 2021.

\bibitem{pacheco2011introduction}
P.S. Pacheco.
\newblock {\em An Introduction to Parallel Programming}.
\newblock Morgan Kaufmann, 2011.

\bibitem{1993HashedOcttreeWarren}
M.S. Warren and J.K. Salmon.
\newblock A parallel hashed oct-tree n-body algorithm.
\newblock In {\em Supercomputing '93:Proceedings of the 1993 ACM/IEEE
  Conference on Supercomputing}, pages 12--21, 1993.

\bibitem{IntelSoA1}
Amanda~K Sharp.
\newblock Memory layout transformations.
\newblock
  \url{https://www.intel.com/content/www/us/en/developer/articles/technical/memory-layout-transformations.html}.
\newblock Accessed: 2021-12-20.

\bibitem{sterling2017high}
T.~Sterling, M.~Brodowicz, and M.~Anderson.
\newblock {\em High Performance Computing: Modern Systems and Practices}.
\newblock Elsevier Science, 2017.

\bibitem{cormen2009introduction}
T.H. Cormen, C.E. Leiserson, R.L. Rivest, and C.~Stein.
\newblock {\em Introduction to Algorithms}.
\newblock MIT Press, 2009.

\bibitem{1990SpatialDataStructuresSamet}
Hanan Samet.
\newblock {\em The Design and Analysis of Spatial Data Structures}.
\newblock Addison-Wesley Longman Publishing Co., Inc., USA, 1990.

\bibitem{ParallelProgramming2013}
Thomas Rauber and Gudula Rünger.
\newblock {\em Parallel Programming}.
\newblock Springer, Berlin, Heidelberg, 2 edition, 2013.

\bibitem{2020LargeFDTD}
Miguel Ruiz-Cabello~N., Maksims Abaļenkovs, Luis~M. Diaz~Angulo, Clemente
  Cobos~Sanchez, Franco Moglie, and Salvador~G. Garcia.
\newblock Performance of parallel fdtd method for shared- and
  distributed-memory architectures: Application tobioelectromagnetics.
\newblock {\em PLOS ONE}, 15(9):1--16, 09 2020.

\end{thebibliography}


\begin{thebibliography}{1}

\bibitem{2019Keyes}
M.~Abduljabbar, M.~A. Farhan, N.~Al-Harthi, R.~Chen, R.~Yokota, H.~Bagci, and
  D.~Keyes.
\newblock Extreme scale fmm-accelerated boundary integral equation solver for
  wave scattering.
\newblock {\em SIAM Journal on Scientific Computing}, 41(3):C245--C268, 2019.

\bibitem{Bauinger2021}
C.~Bauinger and O.~Bruno.
\newblock ``{I}nterpolated {F}actored {G}reen {F}unction'' method for
  accelerated solution of scattering problems.
\newblock {\em Journal of Computational Physics}, 430, 01 2021.

\bibitem{ParallelProgramming2013}
T.~Rauber and G.~Rünger.
\newblock {\em Parallel Programming}.
\newblock Springer, Berlin, Heidelberg, 2 edition, 2013.

\bibitem{sterling2017high}
T.~Sterling, M.~Brodowicz, and M.~Anderson.
\newblock {\em High Performance Computing: Modern Systems and Practices}.
\newblock Elsevier Science, 2017.

\end{thebibliography}

\end{document}


\date{}
\maketitle

\section{Background: computing cluster
  systems} \label{sec:hpcnomenclature} 

The proposed IFGF parallelization strategy is designed for
implementation in modern high performance computing (HPC)
\textit{cluster} systems. The present section reviews relevant
hardware and software concepts and nomenclature utilized throughout
this paper; more detailed descriptions and alternative hardware
designs can be found
e.g. in~\cite{ParallelProgramming2013,sterling2017high}. A modern
computer cluster consists of multiple compute \textit{nodes}. Each
node contains its own memory space, and thus the memory in the cluster
is distributed between the nodes. In particular, access to memory in
other compute nodes requires explicit data communication, which is
typically performed via the \textit{message passing interface}
(MPI)~\cite[Section 8]{sterling2017high}; the performance of
algorithmic implementations for cluster systems can therefore
significantly benefit from careful engineering of MPI-based
inter-process data communications.

Each node typically comprises one or a few \textit{multi-core processors}, each one of which, as the name suggests, contains
multiple computing \textit{cores}. The compute nodes typically are so-called
\textit{shared memory machines} (SMMs), where each core within the
node can access all of the memory in the node. To efficiently make use of more than a single core, certain specialized programming techniques are required, e.g. the Intel Threading Building Blocks (TBB) library, the C++ standard threading model, MPI, or the OpenMP programming interface. Modern compute nodes
usually follow a non-uniform memory access (NUMA) design (in contrast
to uniform memory access (UMA)), where the access times to the shared
memory depend on the locality of the memory with respect to the
multi-core processor accessing it. This design typically results in one or more
NUMA nodes per compute node, where memory access to other
NUMA nodes on the same compute node is usually significantly slower
than access to memory local to the processor.  All of the tests
presented in this paper were conducted on a small cluster consisting
of thirty nodes connected via an InfiniBand interconnect, each one of
which contains four fourteen-core NUMA nodes; additional details
concerning the hardware used are provided in Section~\ref{sec:results}
and Supplementary Materials Sections~\ref{sec:hardware}
and~\ref{sec:pinning}.

On the basis of the functions and synchronization capabilities
provided by MPI, a program can be launched as a set of multiple
\textit{processes} (which are identified in what follows by their
corresponding integer-valued \textit{rank} within the group of all
processes launched by a given program). One of the main roles of the
MPI standard is to allow the programmer to orchestrate the data
communications between the ranks.  Note that, at runtime, an MPI rank can be assigned, or \textit{pinned}, to various
kinds of hardware units, such as e.g. a single core, a NUMA node, a
complete compute node, or various combinations of cores and/or nodes.

\section{Supplementary Numerical Results\label{sec:supp_num_res}}
The present supplementary numerical results section expands on
Section~\ref{sec:results} in various ways. In particular, this section
provides details on the hardware employed and the manner in which it
is used (Sections~\ref{sec:hardware} and~\ref{sec:pinning},
respectively), and, after detailing the method used for error
evaluation (Section~\ref{sec:error}) and the concepts used to evaluate
strong and weak parallel efficiency (Section~\ref{sec:defeff}), it
presents detailed strong parallel scaling results in
Section~\ref{sec:strong} and data regarding the linearithmic
($\mathcal{O}(N\log N)$) scaling tests in
Section~\ref{sec:smlinearithmic}.

\subsection{Compiler and hardware} \label{sec:hardware} The proposed parallel IFGF
program was implemented in C++, and the resulting code was compiled
with the Intel mpiicpc compiler, version 2021.1, and the Intel MPI
library. The following
performance-relevant compiler flags were used: -std=c++20, -O3,
-ffast-math, -qopt-zmm-usage=high, -no-prec-sqrt, -no-prec-div.  All
tests were run on our internal \textit{Wavefield} cluster which
consists of 30 dual-socket nodes.  Each node consists of two Intel
Xeon Platinum 8276 processors with 28 cores per processor, i.e. 56
cores per node, and 384 GB of GDDR4 RAM per node. (The Xeon processors
we use support hyper-threading, but this capability was not exploited
in any of our tests presented in this paper.) The nodes are connected
with HDR Infiniband.

\subsection{Hardware pinning\label{sec:pinning}} As indicated in
Section 3, since each compute node in the \textit{Wavefield} cluster
consists of four NUMA nodes, we run four MPI ranks per node each
pinned to one of these four NUMA nodes through setting the environment
variable ``I\_MPI\_PIN\_DOMAIN=cache3''. The shared memory
parallelization with OpenMP is then used for the parallelization
within each MPI rank, i.e., within a NUMA node. The parallel scaling
within a NUMA node from 1 to 14 cores is investigated below using the
OpenMP specific environment variables ``OMP\_NUM\_THREADS=[1-14]'',
``OMP\_PLACES=cores'', and \\ ``OMP\_PROC\_BIND=true''. The continued
scaling, which is achieved with the MPI parallelization when exceeding
14 cores, is investigated going from one to four MPI ranks (each rank
pinned to one NUMA node in the same compute node), which corresponds
to the MPI scaling on a single, shared-memory node. Finally, the
scaling of the MPI based distributed-memory parallelization is
investigated starting from a single node to 16 nodes, where each node
is fully utilized with four MPI ranks per node and fourteen cores per
rank, as described above. A different hardware pinning used for the
test cases presented in Section~\ref{sec:smlinearithmic} is described
separately in that section.

\subsection{Numerical error estimation\label{sec:error}} The errors
reported in what follows were computed as indicated in
\cite{Bauinger2021}, that is, the relative $L_2$ difference
$\varepsilon_M$ between the full, non-accelerated evaluation of the
field $I(x)$, as stated in \eqref{eq:field1}, and the IFGF-accelerated
evaluation $I_{\text{acc}}(x)$ of \eqref{eq:field1} computed on a
randomly chosen subset of $M$ surface discretization points. More
precisely, the error is given by
\begin{equation} \label{eq:errorcomputation}
 \varepsilon_M = \sqrt{\frac{\sum \limits_{i= 1}^{M} |I(x_{\sigma(i)}) - I_{\text{acc}}(x_{\sigma(i)})|^2}{\sum \limits_{i= 1}^{M} |I(x_{\sigma(i)})|^2}},
\end{equation} 
where $\sigma$ is a random permutation and $x_\ell \in \Gamma_N$
denote the surface discretization points. The method used in
\cite{Bauinger2021} is suitably extended to the present MPI parallel
implementation by using a set of test points $x_\ell$ that contains a
number $M=1000$ of randomly chosen points on each MPI rank. (In
\cite{Bauinger2021} it was shown for sufficiently small examples that
an error evaluation on a subset of 1000 points produces an error
estimation close to the actual error.) More precisely, 1000 surface
discretization points are randomly chosen on each MPI rank from the
distinct set of surface discretization points $\Gamma_{N,\rho}$ each
MPI rank $\rho$ ($1\leq \rho \leq N_r$) is responsible for based on
the distribution introduced in Section
\ref{sec:mpiparallelization}. The final errors are then accumulated
resulting in the overall error estimate
\begin{equation}\label{eq:errorcomputation_2}
  \varepsilon := \varepsilon_M\quad \mbox{ at   $M=1000 \times N_r$ points}.
\end{equation}
As a result, the errors are dependent
on the number $N_r$ of MPI ranks, which is the reason the shown
errors vary slightly as the number of MPI ranks
varies (cf. Tables~\ref{table:strongmpishared} and \ref{table:strongmpidistributed}). All tests were set up in such a way that an error of
approximately $10^{-3}$ is achieved, although, the IFGF method can achieve arbitrarily small errors.

\subsection{Weak and strong parallel efficiency
  concepts} \label{sec:defeff} Let $T(N_c, N)$ denote the time
required by a run of the parallel IFGF algorithm on an $N$-point
discretization $\Gamma_N$, with a given and fixed discretization
scheme, of a given surface $\Gamma$ using $N_c$ cores.  Using this
notation, for a given $N$, the strong parallel efficiency
$E^s_{N_c^0, N_c}$ that results as the number of cores is increased
from $N_c^0$ to $N_c$ is defined as the quotient of the resulting
speedup $S_{N_c^0, N_c}$ to the corresponding ideal speedup value
$S_{N_c^0, N_c}^{\text{ideal}}$:
\begin{align*}
    S_{N_c^0, N_c}^{\text{ideal}} := \frac{N_c}{N_c^0}, \qquad S_{N_c^0, N_c} := \frac{T(N_c^0, \Gamma_N)}{T(N_c, \Gamma_N)}, \qquad E^s_{N_c^0, N_c} := \frac{S_{N_c^0, N_c}}{S_{N_c^0, N_c}^{\text{ideal}}}.
\end{align*}
Note that the implicit dependence on $N$ and $\Gamma_N$ is suppressed
in the speedup and efficiency notations.

The weak parallel efficiency $E^w_{N_c^0, N_c}>0$, in turn, concerns
the computing costs that are observed as the numbers $N_c$ of cores
are increased proportionally to the problem size $N$---effectively
keeping the number of surface discretization points per core
constant---so that as the numbers of cores and discretization points
are simultaneously increased from $N_c^0$ to $N_c$ and from $N^0$ to
$N$, respectively, the relation
\begin{equation}\label{eq:weakscalingrelation}
  N/N^0 = N_c/N_c^0
\end{equation}
is satisfied. Since the weak scaling concerns varying numbers $N$ of
surface discretization points, however, the weak parallel-efficiency
concept must correctly account for the linearithmic theoretical
scaling of the IFGF algorithm. To do this, we consider the computing
time $T(N_c, N)$ required for a run of the algorithm on $N_c$ cores
for an $N$-point discretization of a given, fixed, surface
$\Gamma$. In view of the linearithmic complexity of the algorithm,
perfect weak parallel efficiency would be observed if, for a certain
constant $K$, we had
\[
  T(N_c, N) = \frac{K}{N_c} N \log N.
\]
Thus is to say, under perfect weak parallel scaling, in view
of~\eqref{eq:weakscalingrelation} we would have
\[
  \frac{T(N_c, N)}{T(N^0_c, N^0)} = \frac{N^0_c N\log N }{N_cN^0\log N^0} = \frac{\log N}{\log N^0},
\]
We therefore define the {\em weak parallel efficiency} that results as
the number of cores is increased from $N_c^0$ to to $N_c$ by
\begin{align*}
    E^w_{N_c^0, N_c} :=  \frac{T(N_c^0, N^0) \log{N}}{T(N_c, N) \log{N^0}}.
\end{align*}
Note that $E^w_{N_c^0, N_c} = 1$ corresponds to perfect weak
parallel efficiency, or a weak parallel efficiency of $100\%$.
\begin{table}[H]
\centering
\begin{tabular}{|| c | r | c | c | c | l | r | r ||} 
 \hline
 \multicolumn{1}{|| c |}{$\Gamma$} & \multicolumn{1}{| c |}{\bf $N$} & \multicolumn{1}{|c|}{$d$} & \multicolumn{1}{|c|}{ $N_c$} & \multicolumn{1}{|c|}{$\mathbf \varepsilon$} & \multicolumn{1}{|c|}{ $T$ (s)} & \multicolumn{1}{|c|}{ $E^s_{1, N_c}$ (\%)} & \multicolumn{1}{|c||}{\bf $S_{1, N_c}$} \Tstrut \Bstrut \\ \hline \hline
 \multirow{5}{*}{Sphere} &\multirow{5}{*}{$1,572,864$} & \multirow{5}{*}{$128 \lambda$} & $1$ &  \multirow{5}{*}{$2\times 10^{-3}$} & $2.54\times 10^{3}$ & $100$ & 1.00 \Tstrut \\
  & && $2$ &&  $1.29\times 10^{3}$ & $98$ & 1.95 \Strut \\
  & && $4$ &&  $6.91\times 10^{2}$ & $92$ & 3.67 \Strut \\
  & && $8$ &&  $3.63\times 10^{2}$ & $87$ & 6.98 \Strut \\
  & && $14$ &&  $2.31\times 10^{2}$ & $78$ & 10.98 \Bstrut \\ \hline
  \multirow{5}{*}{\shortstack{Oblate\\Spheroid}} &\multirow{5}{*}{$1,572,864$} & \multirow{5}{*}{$128 \lambda$} & $1$ &  \multirow{5}{*}{$5\times 10^{-4}$} & $9.42\times 10^{2}$ & $100$ & 1.00 \Tstrut \\
  & && $2$ &&  $4.86\times 10^{2}$ & $97$ & 1.94 \Strut \\
  & && $4$ &&  $2.60\times 10^{2}$ & $91$ & 3.62 \Strut \\
  & && $8$ &&  $1.40\times 10^{2}$ & $84$ & 6.69 \Strut \\
  & && $14$ &&  $8.76\times 10^{1}$ & $77$ & 10.75 \Bstrut \\ \hline
  \multirow{5}{*}{\shortstack{Prolate\\Spheroid}} &\multirow{5}{*}{$6,291,456$} & \multirow{5}{*}{$256 \lambda$} & $1$ &  \multirow{5}{*}{$6\times 10^{-4}$} & $1.42\times 10^{3}$ & $100$ & 1.00 \Tstrut \\
  & && $2$ &&  $7.29\times 10^{2}$ & $97$ & 1.95 \Strut \\
  & && $4$ &&  $4.33\times 10^{2}$ & $82$ & 3.28 \Strut \\
  & && $8$ &&  $2.37\times 10^{2}$ & $75$ & 5.99 \Strut \\
  & && $14$ &&  $1.49\times 10^{2}$ & $68$ & 9.49 \Bstrut \\
 \hline
\end{tabular}
\caption{Strong parallel scaling test of the OpenMP IFGF
  implementation from $N_c=1$ to $N_c=14$ cores in a single node for
  three different geometries $\Gamma$. }
\label{table:strongomp128}
\end{table}



\subsection{Strong parallel efficiency tests} \label{sec:strong}

Tables~\ref{table:strongomp128}, \ref{table:strongmpishared}, and
\ref{table:strongmpidistributed} present the strong parallel
efficiencies achieved by the proposed parallel IFGF method under
OpenMP, shared-memory MPI, and distributed-memory MPI parallelization
strategies, respectively, each one for all three test geometries
considered in this section.  In detail, these tables display the main
two strong parallel performance quantifiers, namely the observed
strong parallel efficiency $E^s_{N_c^0, N_c}$ and speedup
$S_{N_c^0, N_c}$, along with the computing times $T$, the obtained
accuracy $\varepsilon$, and details concerning the geometry and the
discretizations. The tables clearly show that, in all cases, the IFGF
parallel efficiencies are essentially independent of the geometry
type. With reference to Section~\ref{sec:pinning} above, the ranges of
the parameter $N_c$ considered in these tables span all of the
available cores in each one of the relevant hardware units used: 14
cores in a single NUMA node, 56 cores (4 NUMA nodes) in a single
compute node, and 16 nodes in the complete cluster (the largest number
of nodes which equals a power of 2 in the cluster used).

\begin{table} [H]
\centering
\begin{tabular}{|| c | r | c | c | c | c | l | r | r ||} 
 \hline
 \multicolumn{1}{|| c|}{$\Gamma$} & \multicolumn{1}{| c |}{\bf $N$} & \multicolumn{1}{|c|}{$d$} & \multicolumn{1}{|c|}{\bf $N_r$} & \multicolumn{1}{|c|}{ $N_c$} & \multicolumn{1}{|c|}{$\mathbf \varepsilon$} & \multicolumn{1}{|c|}{ $T$ (s)} & \multicolumn{1}{|c|}{ $E^s_{14, N_c}$ (\%)} & \multicolumn{1}{|c||}{\bf $S_{14, N_c}$} \Tstrut \Bstrut \\ \hline \hline
 \multirow{3}{*}{Sphere} & \multirow{3}{*}{$1,572,864$} & \multirow{3}{*}{$128 \lambda$} & $1$ & $14$ &  $2\times 10^{-3}$ & $2.31\times 10^{2}$ & $100$ & 1.00 \Tstrut \\
  &&& $2$ &  $28$ & $2\times 10^{-3}$ &  $1.49\times 10^{2}$ & $77$ & 1.54 \Strut \\
  &&& $4$ & $56$ & $2\times 10^{-3}$ &  $7.77\times 10^{1}$ & $74$ & 2.97 \Bstrut \\ \hline
  \multirow{3}{*}{\shortstack{Oblate\\Spheroid}} & \multirow{3}{*}{$1,572,864$} & \multirow{3}{*}{$128 \lambda$} & $1$ & $14$ &  $5\times 10^{-4}$ & $8.76\times 10^{1}$ & $100$ & 1.00 \Tstrut\\
  &&& $2$ &  $28$ & $6\times 10^{-4}$ &  $5.71\times 10^{1}$ & $77$ & 1.53 \Strut \\
  &&& $4$ & $56$ & $6\times 10^{-4}$ &  $2.99\times 10^{1}$ & $73$ & 2.93 \Bstrut \\ \hline
  \multirow{3}{*}{\shortstack{Prolate\\Spheroid}} & \multirow{3}{*}{$6,291,456$} & \multirow{3}{*}{$256 \lambda$} & $1$ & $14$ &  $6\times 10^{-4}$ & $1.49\times 10^{2}$ & $100$ & 1.00 \Tstrut\\
  &&& $2$ &  $28$ & $6\times 10^{-4}$ &  $9.10\times 10^{1}$ & $82$ & 1.65 \Strut \\
  &&& $4$ & $56$ & $5\times 10^{-4}$ &  $4.97\times 10^{1}$ & $75$ & 3.01 \Bstrut \\
 \hline
\end{tabular}
\caption{Strong parallel scaling test of the shared-memory MPI
  implementation on a single node, transitioning from $N_c=14$ cores
  to $N_c=56$ (all cores available in one compute node) by increasing
  the number $N_r$ of MPI ranks from 1 to 4, for three different
  geometries $\Gamma$. }
\label{table:strongmpishared}
\end{table}

The largest efficiency deficit observed as a result of a
hardware-doubling transition is the decrease by a full $23\%$ (from
100\% to 77\%) shown in Table~\ref{table:strongmpishared}, which
results from the transition from one to two MPI ranks (that is, from
one to two 14-core NUMA nodes). We argue that this deficit, which
takes place precisely as an MPI communication between NUMA nodes is
first introduced, is not a sole reflection of the character of the
algorithm in presence of the MPI interface, since such large deficits
are not observed in any other MPI related hardware-doubling
transitions reported in the various tables. As potential additional
contributing elements to this deficit we mention notably, MPI overhead
(which would only be incurred in the first doubling transition but not
in subsequent doubling transitions, in view of the decreasing number
of pairwise communications incurred by the algorithm under a doubling
transition in a strong scaling test, as indicated by the theoretical
discussion in Section~\ref{sec:complexity}), and the Intel Turbo Boost
Technology inherent in the processors used---which achieve maximum
turbo frequencies when running under lower loads, and which, when
concurrently using larger numbers of cores in a single node, cease to
operate.

A variety of other data is presented in these
tables. Tables~\ref{table:strongomp128}
and~\ref{table:strongmpishared}, which demonstrate the strong scaling
within a single NUMA node, and among all four NUMA nodes within a
compute node, are included for completeness, but as discussed below,
we attach far greater significance to
Table~\ref{table:strongmpidistributed}, which demonstrates the scaling
of the method under the one hardware element that can truly be
increased without bounds, namely, the number of compute nodes.  In
this table, geometries twice as large than those used for the previous
two tables are considered (to reasonably increase the minimum
computing times), and the hardware is scaled from one compute node to
sixteen compute nodes. Per the description in
Section~\ref{sec:pinning}, each node is assigned four MPI ranks, each
one of which is pinned to one of the four NUMA nodes present in the
compute node. Overall, a strong scaling efficiency of over 60\% can be
observed in all cases, with the results of the sphere test case even
above 70\% owing to the symmetry of the geometry and the resulting
increased load-balance and minimized communication between ranks. The
loss of efficiency can be attributed the load-imbalance induced by our
data partitioning strategy, the communication between ranks, and the
parallelization overhead introduced by MPI and OpenMP.

\begin{table} [H]
\centering
\begin{tabular}{|| c | r | c | c | l | l | r | r | r ||} 
 \hline
 \multicolumn{1}{|| c|}{$\Gamma$} & \multicolumn{1}{| c |}{\bf $N$} & \multicolumn{1}{|c|}{$d$}  & \multicolumn{1}{|c|}{ $N_c$} & \multicolumn{1}{|c|}{$\mathbf \varepsilon$} & \multicolumn{1}{|c|}{$T$ (s)} & \multicolumn{1}{|c|}{\bf $E^s_{56, N_c}$ (\%)} & \multicolumn{1}{|c|}{\bf $S_{56, N_c}$} & \multicolumn{1}{|c||}{\bf $E^s_{\frac{N_c}{2}, N_c}$}  \Tstrut \Bstrut \\ \hline \hline
 \multirow{5}{*}{Sphere} & \multirow{5}{*}{$6,291,456$} & \multirow{5}{*}{$256 \lambda$} & $56$ &  $2 \times 10^{-3}$ & $3.55\times 10^{2}$ & $100$ & 1.00 & -  \Tstrut \\
  &&& $112$ & $2\times 10^{-3}$ &  $1.80\times 10^{2}$ & $99$ & 1.97 & $99$ \Strut \\
  &&& $224$ & $2\times 10^{-3}$ &  $9.78\times 10^{1}$ & $91$ & 3.64 & $92$ \Strut \\
  &&& $448$ & $2\times 10^{-3}$ &  $5.52\times 10^{1}$ & $81$ & 6.44 & $89$ \Strut \\
  &&& $896$ & $2\times 10^{-3}$ & $3.12\times 10^{1}$ & $71$ & 11.40 & $89$ \Bstrut \\ \hline
  \multirow{5}{*}{\shortstack{Oblate\\Spheroid}} & \multirow{5}{*}{$6,291,456$} & \multirow{5}{*}{$256 \lambda$} & $56$ &  $6\times 10^{-4}$ & $1.34\times 10^{2}$ & $100$ & 1.00 & - \Tstrut \\
  &&& $112$ & $6\times 10^{-4}$ &  $7.41\times 10^{1}$ & $91$ & 1.81 & $91$ \Strut \\
  &&& $224$ & $6\times 10^{-4}$ &  $4.17\times 10^{1}$ & $80$ & 3.22 & $89$ \Strut \\
  &&& $448$ & $6\times 10^{-4}$ &  $2.38\times 10^{1}$ & $70$ & 5.64 & $88$ \Strut \\
  &&& $896$ & $6\times 10^{-4}$ & $1.40\times 10^{1}$ & $60$ & 9.56 & $85$ \Bstrut \\ \hline
  \multirow{5}{*}{\shortstack{Prolate\\Spheroid}} & \multirow{5}{*}{$25,165,824$} & \multirow{5}{*}{$512 \lambda$} & $56$ &  $4\times 10^{-4}$ & $2.23\times 10^{2}$ & $100$ & 1.00 & - \Tstrut \\
  &&& $112$ & $5\times 10^{-4}$ &  $1.22\times 10^{2}$ & $92$ & 1.83 & $92$ \Strut \\
  &&& $224$ & $6\times 10^{-4}$ &  $6.83\times 10^{1}$ & $82$ & 3.28 & $89$ \Strut \\
  &&& $448$ & $6\times 10^{-4}$ &  $3.74\times 10^{1}$ & $75$ & 5.97 & $91$ \Strut \\
  &&& $896$ & $6\times 10^{-4}$ & $2.23\times 10^{1}$ & $63$ & 10.01 & $84$ \Bstrut \\
 \hline
\end{tabular}
\caption{Strong parallel scaling test of the distributed-memory MPI
  implementation from $N_c=56$ to $N_c=896$ cores (1 to 16 compute
  nodes) with 4 MPI ranks per node for three different geometries
  $\Gamma$. }
\label{table:strongmpidistributed}
\end{table} 
The most important quality illustrated in these tables is the IFGF's
efficiency performance under strong-scaling hardware-doubling
transitions demonstrated in the last column of
Table~\ref{table:strongmpidistributed}. This performance, which
mirrors the corresponding weak-scaling performance presented in the
last column of Table~\ref{table:weakmpi}, shows that, as in the weak
scaling case, under the assumption that the displayed trend is
maintained for large numbers of nodes (below the obvious limit
imposed by the fixed problem size), the parallel IFGF method for a
fixed problem can be efficiently run in large numbers of computing
cores---with an efficiency factor no worse than a constant
$\approx 80\%$ as the hardware sizes are doubled from a given point of
reference.

The character of the IFGF algorithm under weak- and strong-scaling
hardware-doubling tests, as discussed above in this section and in
Section~\ref{sec:complexity}, would ensure that, provided the
demonstrated trends are maintained (as is expected in view of the
discussion in the first paragraph of
Section~\ref{sec:mpiparallelization}), the method can be executed
successfully in very large hardware infrastructures.

\subsection{Linearithmic scaling}\label{sec:smlinearithmic} Table~\ref{table:linearithmic} presents numerical data related to the test case considered
in Figure~\ref{fig:linearithmic}, which concerns the linearithmic
scaling of the parallel IFGF method on a fixed and large number of
cores.  As indicated in Section~\ref{sec:results}, the data in the
figure and table were generated by pinning a single MPI rank to each
compute node, each of which spawns 56 OpenMP threads on each one of
the 30 nodes available in the computer cluster we use, with parameters
resulting in an IFGF error $\varepsilon\approx 1.5\times 10^{-2}$
(cf. equation~\eqref{eq:errorcomputation_2}). The observed complexity
slightly outperforms the postulated $\mathcal{O}(N \log N)$ estimate
within this range of values of $N$. Note, in particular, the last
column of Table~\ref{table:linearithmic} which suggests rapid
convergence to exact linearithmic complexity, with a well defined
proportionality constant, as $N$ grows.

\begin{table} [ht]
\centering
\begin{tabular}{|| c | r | r | c | l | r | r | r ||} 
 \hline
\multicolumn{1}{|| c|}{$\Gamma$} & \multicolumn{1}{| c |}{\bf $N$} & \multicolumn{1}{|c|}{$d$} & \multicolumn{1}{|c|}{\bf $N_c$} & \multicolumn{1}{|c|}{$\mathbf \varepsilon$} & \multicolumn{1}{|c|}{$T$ (s)} & \multicolumn{1}{|c|}{Mem/rank} & \multicolumn{1}{|c|}{$T/(N \log N)$}  \Tstrut \Bstrut \\ \hline \hline
\multirow{4}{*}{\shortstack{Prolate\\Spheroid}} & $6,291,456$ &$512 \lambda$ & \multirow{4}{*}{$1,680$} &  \multirow{4}{*}{$1.5\times10^{-2}$} & $4.09\times 10^{0}$ & $0.50$ GB & $9.56 \times 10^{-8}$ \Tstrut \\
 & $25,165,824$ &$1,024 \lambda$& & & $1.64\times 10^{1}$ & $1.89$  GB & $8.83 \times 10^{-8}$ \Strut \\
  &$100,663,296$ & $2,048 \lambda$& & &  $6.71\times 10^{1}$ & $7.42$  GB & $8.33 \times 10^{-8}$ \Strut \\
  &$402,653,184$ & $4,096 \lambda$& & &  $2.90\times 10^{2}$ & $29.73$  GB & $8.37 \times 10^{-8}$ \Bstrut \\
 \hline
\end{tabular}
\caption{Preservation of the linearithmic IFGF scaling in the parallel
  context.  One MPI rank per node and 56 OpenMP threads per MPI rank
  on 30 compute nodes for a prolate spheroid geometry were used for
  this test. The peak memory per MPI rank used by the IFGF method
  (excluding the memory required to store the initial geometry) is
  listed in the next-to-last column. The value in the last column
  suggests convergence to exact $\mathcal{O}(N \log N)$ scaling.}
\label{table:linearithmic}
\end{table}

\subsection{Large sphere tests}\label{sec:largespheretests}
Table~\ref{table:largesphere} illustrates the performance of the IFGF
method in terms of computing time and memory requirements for several
large-sphere configurations, all of them run on our full 30 node,
1,680 core cluster. In particular, Table~\ref{table:largesphere} shows
that, as mentioned in Section~\ref{sec:introduction}, on the basis of
approximately 4 TB of memory, a sphere 1,389 $\lambda$ in diameter
(which, as described below, corresponds to the case $f=238.086$KHz
considered in~\cite{2019Keyes}) with 2.12 billion DOF was run in a
computing time of 2,380 seconds with a 0.5$\%$ near-field error (which
may be compared to the only error indicator reported in \cite[Table
2]{2019Keyes} for this test case, which amounts to 20$\%$, and
compared to the 3$\%$ near-field solution error reported in the same
table of that paper for significantly smaller problems). The computing
time reported for this test case is a factor of approximately 46 times
longer than the time reported in~\cite{2019Keyes}, for the same number
of DOFs and sphere size, on a computer 78 times larger (containing
131,072 cores) and on the basis of an unspecified amount of memory.

\begin{table} [H]
\centering
\begin{tabular}{|| c | r | r | c | l | r | r ||} 
 \hline
\multicolumn{1}{|| c|}{$\Gamma$} & \multicolumn{1}{| c |}{\bf $N$} & \multicolumn{1}{|c|}{$d$} & \multicolumn{1}{|c|}{\bf $N_c$} & \multicolumn{1}{|c|}{$\mathbf \varepsilon$} & \multicolumn{1}{|c|}{$T$ (s)} & \multicolumn{1}{|c|}{Mem/rank}  \Tstrut \Bstrut \\ \hline \hline
  \multirow{5}{*}{\shortstack{Sphere}}
& $1,610,612,736$ &$1,389 \lambda$&\multirow{5}{*}{$1,680$} & $4\times10^{-3}$ & $3.59\times 10^{3}$ & $125.96$ GB \Tstrut \\
& $1,610,612,736$ &$2,048 \lambda$ & &  $7\times10^{-2}$ & $2.84\times 10^{3}$ & $105.91$ GB \Strut \\
&$1,944,000,000$ & $1,389 \lambda$& & $1\times10^{-1}$ &  $1.01\times 10^{3}$ & $42.45$  GB \Strut \\
&$1,944,000,000$ & $1,389 \lambda$& & $5\times10^{-3}$ &  $2.34\times 10^{3}$ & $133.44$  GB \Strut \\
&$2,120,640,000$ & $1,389 \lambda$& & $5\times10^{-3}$ &  $2.38\times 10^{3}$ & $134.63$  GB \Bstrut \\
 \hline
\end{tabular}
\caption{Large sphere test cases run on thirty 56-core compute nodes
  (for a total of $1,680$ cores), utilizing thirty MPI ranks. The
  sphere of acoustic size $1,389 \lambda$ in this table coincides with
  largest sphere test case considered in~\cite{2019Keyes}. }
\label{table:largesphere}
\end{table}

As indicated above, the sphere 1,389 $\lambda$ in diameter in this
table coincides with largest sphere test case considered
in~\cite{2019Keyes}, which is cited in Table~2 in that reference as a
sphere of two-meters in diameter illuminated at the frequency of
$f= 238.086$ KHz under the speed of sound $c = 343m/s$ assumed in that
reference. (The number of wavelengths in this two-meter sphere test
case equals $2$
meters$/\lambda = 2 f/c = 2\cdot 238,086/343 \approx 1,388.26 <
1,389$.)  The largest discretization presented in the present table
for this sphere test-case ($2,120,640,000$ discretization points, a
limit induced by the largest number representable by a signed 32-bit
integer assumed in our geometry-generation code, which will be avoided
in subsequent code implementations by switching to 64-bit integers),
is slightly smaller than the $2,300,067,840$ discretization considered
in~\cite{2019Keyes} under a 131,072 core run.

Other test cases listed in Table~\ref{table:largesphere} include an
example for a much larger sphere, $2,048\lambda$ in diameter, as well
as other $1,389 \lambda$ test cases for various accuracies and
discretization sizes---and, in all cases, on the basis of memory
consumptions ranging between $\approx 1.2$TB and $\approx 4$TB.

\bibliographystyle{abbrv}
\bibliography{IFGF}